\documentclass{article}%
\usepackage{fullpage}
\usepackage{amsfonts}
\usepackage{amssymb}
\usepackage{amsthm}
\usepackage{amsmath}
\usepackage{float}
\usepackage{arydshln}
\usepackage{authblk}
\usepackage{mleftright}
\usepackage{graphicx}
\usepackage{tikz}
\usetikzlibrary{decorations.pathreplacing}
\usepackage{xcolor}%
\providecommand{\U}[1]{\protect\rule{.1in}{.1in}}
\newtheorem{theorem}{Theorem}
\newtheorem{lemma}[theorem]{Lemma}

\newtheorem{corollary}[theorem]{Corollary}

\newtheorem{definition}[theorem]{Definition}
\newtheorem{notation}[theorem]{Notation}
\newtheorem{algorithm}[theorem]{Algorithm}
\newtheorem{example}[theorem]{Example}
\newtheorem{remark}[theorem]{Remark}
\newtheorem{problem}[theorem]{Problem}
\allowdisplaybreaks
\DeclareMathOperator{\sseq}{sseq}
\DeclareMathOperator{\cl}{cl}
\DeclareMathOperator{\one}{\delta}
\DeclareMathOperator{\diag}{diag}
\newcommand{\As}{A_{\rm{sig}}}
\newcommand{\Cs}{C_{\rm{sig}}}
\newcommand{\EC}{{\rm{EC}}}
\newcommand{\GEC}{{\rm{GEC}}}

\newcommand{\sign}{\operatorname*{sign}}

\author[1]{Hoon Hong\thanks{hong@ncsu.edu}}
\author[1]{Daniel Profili\thanks{daprofil@ncsu.edu, corresponding author}}
\author[2]{J. Rafael Sendra\thanks{jrafael.sendra@cunef.edu}}

\affil[1]{Department of Mathematics, North Carolina State University, USA}

\affil[2]{Department of Mathematics, CUNEF Universidad, Spain}

\begin{document}

\title{Conditions for eigenvalue configurations \\of two real symmetric matrices \\ (Signature approach) \thanks{A special case of the result here, without proofs, was presented as a poster at ISSAC 2024 with an abstract in \cite{issac} }}
\maketitle

\begin{abstract} 
For two real symmetric matrices, their eigenvalue configuration is the
relative arrangement of their eigenvalues on the real line.
We consider the following problem: given two parametric real symmetric matrices and an eigenvalue configuration, find a
simple condition on the parameters such that the two matrices have
the given eigenvalue configuration.
In this paper, we develop theory and give an algorithm for this problem, { allowing the case where the matrices have multiple or shared eigenvalues.}
The output of the algorithm is a condition written in terms of the signatures of certain related symmetric matrices.
We give a concrete application to parametric reachability of rank-$r$ updates of real symmetric matrices.
\end{abstract}

\section{Introduction} 
For two real symmetric matrices $F$~and~$G$, their eigenvalue configuration is the
relative arrangement of their eigenvalues on the real line.
{ In the simplest case, consider two real symmetric matrices $F$ and $G$ with respective simple eigenvalues $\alpha_1 < \alpha_2$ and $\beta_1 < \beta_2 < \beta_3$ where their eigenvalues are ordered as in the following picture.
  \[
  \begin{tikzpicture}
    \draw (0, 0) -- (7.5, 0);
    \fill [red] (1,0) circle(3pt)
    node[label=$\alpha_1$]{};
    \fill [red] (5,0) circle(3pt)
    node[label=$\alpha_2$]{};
    \fill [blue] (3,0) circle(3pt)
    node[label=$\beta_1$]{};
    \fill [blue] (4,0) circle(3pt)
    node[label=$\beta_2$]{};
    \fill [blue] (6,0) circle(3pt)
    node[label=$\beta_3$]{};
  \end{tikzpicture}
    \]
    When the eigenvalues are simple and distinct, we define eigenvalue configuration as the vector of numbers of eigenvalues of $G$ ($\beta$'s) which lie in the intervals between each pair of consecutive eigenvalues of $F$ ($\alpha$'s). In the example above, there are two eigenvalues of $G$ ($\beta_1$ and $\beta_2$) in the interval $(\alpha_1, \alpha_2)$, and one eigenvalue of $G$ ($\beta_3$) in the interval $(\alpha_2, \infty)$, so the eigenvalue configuration is the vector $(2,1)$.
    In the case where the matrices have multiple eigenvalues, we simply count with multiplicity, and in the case where eigenvalues of $F$ and $G$ coincide, the eigenvalues of $G$ are shared between neighboring intervals in a natural way similar to techniques which appear in complex root counting \cite{mariealgebraic}.
In the next section, we give a precise definition of eigenvalue configuration.
} In this paper, we tackle the following problem: 
Given two parametric real symmetric matrices and an eigenvalue  configuration, find  a quantifier-free necessary and sufficient condition on the parameters so that the eigenvalue configuration of the matrices is
the given one. 

For some historical context, a fundamental problem in computational algebra and geometry is to 
determine
quantifier-free necessary and sufficient conditions on the coefficients of a
polynomial such that its roots lie in a given subset of the complex plane. This
problem appears in a wide breadth of fields, including the theory of stable polynomials \cite{michelen2024central},
combinatorics \cite{codenotti2022combinatorics}, graph theory
\cite{gurvits2023trees}, and functional analysis \cite{aleman2023cyclicity}, among others.
The eigenvalue configuration problem is a specific case of this problem.

The eigenvalue configuration problem generalizes Descartes' rule of signs, which is a fundamental and widely used (see e.g. \cite{chattopadhyay2023near}, \cite{haija2021resource}, \cite{Huang_2021}) tool in computational
real algebraic geometry.  
Recall that Descartes' rule of signs states that, for a real polynomial
$g$, the number of
positive real roots of $g$ (counted with multiplicity) is bounded above by the sign variation count of the coefficients of
$g$; i.e., the number of times consecutive coefficients change sign, ignoring
zeros. It is especially useful when the polynomial $g$ has only real roots since the number of positive real roots of $g$
is exactly the sign variation count of the coefficients.  Viewed as an eigenvalue configuration problem, Descartes' rule of signs
can be thought of as determining the eigenvalue configuration of the $1
\times1$ matrix $F = [0]$ and a symmetric matrix~$G$ whose characteristic
polynomial is $g$. 
{ Therefore, the eigenvalue configuration problem extends Descartes' rule of signs by generalizing the \textit{positive root counting} of Descartes' rule of signs into a \textit{configuration problem}, where we discriminate different arrangements of the roots of two polynomials instead of just one. To simplify the presentation, in this paper we reframe the problem slightly in terms of characteristic polynomials of real symmetric (or Hermitian) matrices.}
Given the wide applicability of Descartes' rule of signs, we expect that a generalization may have even more applications. For
an example, one could use this generalization in investigating the impact on the eigenvalues under low rank updates~\cite{benaych2011eigenvalues,hill1992refined}. 

The main difficulty of the eigenvalue configuration problem comes from the fact that it is not practically solvable using existing { general} methods.
It is true that the eigenvalue configuration problem can be, in principle, solved using general quantifier elimination algorithms 
(see e.g. \cite{McCallum:84,
BKR1986, Grigorev:88,Weispfenning:95,
Gonzalez_Lombardi_Recio_Roy:89, Hong:90a, Hong:90b, Collins_Hong:91, Hong:92a,
Renegar:92a,Canny:93a,Loos_Weispfenning:93, Weispfenning:94b,
Gonzalez-Vega:96, DBLP:journals/jsc/McCallum99, DBLP:conf/issac/McCallumB09,
DBLP:journals/jsc/Brown01a, DBLP:journals/jsc/Brown01, Strzebonski06,
DBLP:conf/issac/ChenMXY09, DBLP:journals/jsc/Brown12, Hong_Safey:2012}), which are readily available as free software tools.
Thus, one might wonder whether we could just run those software tools on the
above quantified formula. However, there are three major difficulties: (1) It
requires that we run the software on each $m$ and $n$ from scratch (where $m$ and $n$ are the sizes of the two input matrices).
(2) Each run (on each $m$ and $n$) quickly becomes very time-consuming as
$m$ and $n$ grow, making this approach practically infeasible and (3) the
output formula looks chaotic, making it practically incomprehensible.
As a consequence of these limitations, { we need to develop a specialized approach which exploits the nature of the eigenvalue configuration problem.}
{
  The nontriviality of the eigenvalue configuration problem comes from the fact that, in general, one cannot symbolically represent the eigenvalues of a matrix using the language of semialgebraic sets: addition, multiplication, equalities, and inequalities. 
In addition, since our goal is to
parametrically characterize \textit{all} possible matrices whose eigenvalues
are arranged in a certain way, numerical approaches cannot be used. It is useful to have such a semialgebraic representation for many reasons, including ease of reasoning and manipulation as well as use as inputs for other concrete applications or computer algebra systems.
  }

The main contribution of this paper is to provide an efficient and structured solution to the eigenvalue configuration problem. We accomplish this by producing a set { of $2^m$ matrices, each $n \times n$}, whose entries are polynomials in the parameters of the input matrices. The signatures of these new matrices determine the eigenvalue configuration via a linear map
(depending only on the sizes of $F$~and~$G$ and not on their
parameters) which takes a column vector of integer signatures and maps to an eigenvalue configuration. {This paper crucially uses the generalized root counting theorem given by Ben-Or, Kozen, and Reif \cite{BKR1986}, but our main contribution is neither a modification nor an extension of that result. We develop new combinatorial machinery for studying eigenvalue configurations specifically, and BKR's theorem is used only as a black box (for details, see Lemma \ref{lem:elim_count_N} and its proof).} In our related work \cite{symmetric}, we approach the same eigenvalue configuration problem via the Fundamental Theorem of Symmetric Polynomials.



As discussed, our contribution can be seen as one possible way of generalizing
Descartes' rule of signs to { a configuration problem on the roots of two polynomials}. One could also generalize Descartes' rule of signs in the other direction; i.e., by increasing the number of variables rather than the number of polynomials. {Descartes' rule of signs has been generalized to multivariate systems (\cite{bihan2017descartes, itenberg1996multivariate}) as well as single multivariate polynomials (\cite{telek2024geometry}). These differ from our result because they treat either solutions to polynomial systems or a single multivariate polynomial, rather than a configuration of roots of multiple polynomials.} To our knowledge, there has not yet
been significant work on generalizing Descartes' rule of signs for { root configurations of} multiple univariate
polynomials, { which is the problem we tackle in this paper.}


{ 
  The structure of the paper is as follows. First, in Section \ref{sec:ec} we discuss and precisely define the eigenvalue configuration of two real symmetric matrices.
  In Section \ref{sec:problem}, we state the problem. 
  In Section \ref{sec:main}, we state our main
  theorem (Theorem \ref{thm:main}). The proof (Section \ref{sec:proof}) is
  divided into three parts. First, in Section~\ref{sec:strong}, we prove Theorem
  \ref{thm:main} for a special case. Following
  that, in Section \ref{sec:generic}, we generalize the result slightly. In Section~\ref{sec:arbitrary}, we prove our theorem for arbitrary real symmetric matrices. {To conclude, in Section~\ref{sec:update} we briefly discuss a potential concrete application of our main theorem to the reachability of rank-$r$ updates.}
}

\section{Problem}

{ 
  \subsection{Defining eigenvalue configuration}
  \label{sec:ec}

  In this section, we will precisely define what ``eigenvalue configuration'' means. 
Let $F \in\mathbb{R}^{m \times m}$ and $G \in\mathbb{R}^{n \times n}$ be real
symmetric matrices. {Note that ``real symmetric'' could be substituted with ``Hermitian'' throughout the paper with no effect on the result, as we only require that all eigenvalues are real. For the simplicity of presentation, we will focus on the real case.}
First, let us consider the simplest case when $F$~and~$G$ have only simple eigenvalues and do not share any eigenvalues.
Since $F$~and~$G$ are real symmetric, all their eigenvalues are real. This means we can meaningfully consider the ordering of both matrices' eigenvalues on the real line.

Consider the following example.
Suppose $F$~and~$G$ are $3 \times 3$ real symmetric matrices such that their eigenvalues are arranged as in the following diagram, 
where the red points denote eigenvalues of $F$ and blue points denote eigenvalues of $G$.
\[
  \begin{tikzpicture}
    \draw (-2, 0) -- (5, 0); \fill [red] (0,0) circle(3pt) node[]{};
    \fill [red] (1,0) circle(3pt) node[]{}; \fill [red] (3,0)
    circle(3pt) node[]{}; \fill [blue] (-1,0) circle(3pt) node[]{};
    \fill [blue] (2,0) circle(3pt) node[]{}; \fill [blue] (4,0)
    circle(3pt) node[]{};
  \end{tikzpicture}
  \]
  
Our goal is to define a notion of eigenvalue configuration which encapsulates the information depicted by the diagram above. 
One natural way to do this is to count the number of blue points between each pair of red points.
In this example, there are 3 red points, and so they partition the real line into 4 open intervals. Then, we count the number of blue points in each interval, as in
the following diagram.
\[
  \begin{tikzpicture}
    \draw (-2, 0) -- (5, 0); \fill [red] (0,0) circle(3pt) node[]{};
    \fill [red] (1,0) circle(3pt) node[]{}; \fill [red] (3,0)
    circle(3pt) node[]{}; \fill [blue] (-1,0) circle(3pt) node[]{};
    \fill [blue] (2,0) circle(3pt) node[]{}; \fill [blue] (4,0)
    circle(3pt) node[]{}; \draw
    [decorate,decoration={brace,amplitude=10pt},xshift=0pt,yshift=-10pt]
    (0,0) -- (-2,0) node [black,midway,xshift=0cm,yshift=-0.75cm]
    {$\# \tikz\draw[blue,fill=blue] (0,0) circle (.5ex); = 1$}; \draw
    [decorate,decoration={brace,amplitude=10pt},xshift=0pt,yshift=-10pt]
    (1,0) -- (0,0) node [black,midway,xshift=0cm,yshift=-0.75cm]
    {$\# \tikz\draw[blue,fill=blue] (0,0) circle (.5ex); = 0$}; \draw
    [decorate,decoration={brace,amplitude=10pt},xshift=0pt,yshift=-10pt]
    (3,0) -- (1,0) node [black,midway,xshift=0cm,yshift=-0.75cm]
    {$\# \tikz\draw[blue,fill=blue] (0,0) circle (.5ex); = 1$}; \draw
    [decorate,decoration={brace,amplitude=10pt},xshift=0pt,yshift=-10pt]
    (5,0) -- (3,0) node [black,midway,xshift=0cm,yshift=-0.75cm]
    {$\# \tikz\draw[blue,fill=blue] (0,0) circle (.5ex); = 1$};
  \end{tikzpicture}
  \]
  
The number of blue points in each of these intervals is, respectively, $(1, 0, 1, 1)$.
Note that the sum of the entries of this vector is exactly the number of eigenvalues of $G$.
Since all roots of $G$ are real, we can save some space by omitting one of the counts; by arbitrary choice, we can omit the first count (i.e. the number of blue points which lie to the left of the first red point), which could later be reobtained by adding up the remaining entries and subtracting from the size of $G$.
Thus, in this case, we will encode the eigenvalue configuration as $(0,1,1)$.

Note that this definition is inherently biased. We could just as easily have instead counted the number of red points between each pair of blue points. This choice is mostly arbitrary; however, in many potential applications (some of which we discuss in the next section), we are interested in the eigenvalue configuration of a matrix $F$ and another matrix derived by $F$; for example, one might consider a perturbation of the entries of $F$ to study how the eigenvalues change. In these contexts, it is natural to bias our definition in favor of~$F$ as the ``base'' matrix.

\bigskip
\noindent Now, we will make this definition precise.
First, we introduce some notation.

\begin{notation}
\label{basic_notation}\ \ 
\end{notation}

\begin{enumerate}
\item Let $F = [a_{ij}] \in\mathbb{R}^{m \times m}$ and $G = [b_{ij}] \in\mathbb{R}^{n \times n}$ be
real symmetric matrices.

\item Let $\alpha=\left(  \alpha_{1},\ldots,\alpha_{m}\right)  $ be the
eigenvalues of$\ F$ and let $\beta=\left(  \beta_{1},\ldots,\beta_{n}\right)  $ \ be the eigenvalues
of $G$.

Since $F$~and~$G$ are real symmetric, all their eigenvalues are real. Thus without losing
generality, let us index the eigenvalues so that
$\alpha_{1}\le\alpha_{2}\le\cdots\le\alpha_{m}$ and $\beta_{1}  \le\beta_{2}\le\cdots\le\beta_{n}$.

\item Let $A_{t}$ denote the { interval} $\{x \in \mathbb{R}: \alpha_{t} < x < \alpha_{t+1}\}$ for $t=1,\ldots,m$, where
$\alpha_{m+1}=\infty$. (Note that this interval could be empty if $\alpha_t = \alpha_{t+1}$.)
\end{enumerate}
We have already addressed the case where $F$~and~$G$ have simple and distinct eigenvalues; in that case, the eigenvalue configuration is the vector $c = (c_1, \dots, c_m)$ where $c_t = \# \{i : \beta_i \in A_t \}$.

Next, let us generalize slightly to the case where $F$~and~$G$ still do not share eigenvalues, but multiple eigenvalues are allowed.
We will refer to this case often in the rest of the paper, so we will give a special name.

\begin{definition}
  [Generic]\label{def:generic} We say that the pair $(F,G)$ of real symmetric matrices is \textbf{generic} if
  $F$~and~$G$ do not share eigenvalues.

  \medskip

  \noindent Note that this definition describes a property of {\rm pairs} of matrices.
  We may sometimes be more casual with this phrasing and simply say that ``$F$ and $G$ are generic matrices.''
\end{definition}
\noindent The motivation for this terminology is the fact that if $(F,G)$ is a generic pair, then sufficiently small perturbations of the eigenvalues of $F$~and~$G$ preserve the relative arrangement of their eigenvalues.

In the case of generic matrices which \textbf{do} contain multiple (repeated) eigenvalues, we can use the previous definition with one minor caveat.
For instance, suppose $F$~and~$G$ are such that their eigenvalues are arranged as in
the following diagram.

\[
  \begin{tikzpicture}
    \draw (-1, 0) -- (5, 0); \fill [red] (0,0) circle(3pt) node[]{};
    \fill [red] (0,0.5) circle(3pt) node[]{}; \fill [red] (3,0)
    circle(3pt) node[]{}; \fill [blue] (2,0) circle(3pt) node[]{};
    \fill [blue] (2,0.5) circle(3pt) node[]{}; \fill [blue] (4,0)
    circle(3pt) node[]{};
  \end{tikzpicture}
  \]
  
In this case, the only extra consideration we need to make is that the interval between the first two red points is the empty set, so the corresponding entry in the eigenvalue configuration vector will be zero.
This leads to the following definition for the eigenvalue configuration of generic matrices.

\begin{definition}
[Eigenvalue configuration for generic matrices]\label{def:GEC}The eigenvalue configuration of generic pairs of matrices $F$ and~$G$, written as $\GEC\left(  F,G\right)$, is the
tuple $$c=(c_{1},\dots,c_{m})$$ where $$c_{t}=\#\{i:\beta_{i}\in A_{t}\}.$$
\end{definition}

\begin{example}[Eigenvalue configuration of generic matrices]
Let $F \in\mathbb{R}^{5 \times5}$ and $G \in\mathbb{R}^{7 \times7} $ be symmetric
matrices such that their corresponding eigenvalues are
\[
\alpha=(0,0,3,9,12),\qquad\beta=(-1,2,4,4,8,10,10).
\]
Note that the pair $(F,G)$ is generic since $F$ and $G$ do not share any eigenvalues.
The eigenvalues are arranged on the real line as in the following diagram
\[
  \begin{tikzpicture}
    \draw (-2, 0) -- (13, 0); \fill [red] (0,0) circle(3pt)
    node[label=$\alpha_1$]{}; \fill [red] (0,0.8) circle(3pt)
    node[label=$\alpha_2$]{}; \fill [red] (3,0) circle(3pt)
    node[label=$\alpha_3$]{}; \fill [red] (9,0) circle(3pt)
    node[label=$\alpha_4$]{}; \fill [red] (12,0) circle(3pt)
    node[label=$\alpha_5$]{}; \fill [blue] (-1,0) circle(3pt)
    node[label=$\beta_1$ ]{}; \fill [blue] (2,0) circle(3pt)
    node[label=$\beta_2$ ]{}; \fill [blue] (4,0) circle(3pt)
    node[label=$\beta_3$ ]{}; \fill [blue] (4,0.8) circle(3pt)
    node[label=$\beta_4$ ]{}; \fill [blue] (8,0) circle(3pt)
    node[label=$\beta_5$]{}; \fill [blue] (10,0) circle(3pt)
    node[label=$\beta_6$]{}; \fill [blue] (10,0.8) circle(3pt)
    node[label=$\beta_7$]{};
  \end{tikzpicture}
  \]
  

\medskip

\noindent Then
\[%
\begin{array}
[c]{ll}%
A_{1}=(\alpha_{1},\alpha_{2}) & \\
A_{2}=(\alpha_{2},\alpha_{3}) & \ni\beta_{2}\\
A_{3}=(\alpha_{3},\alpha_{4}) & \ni\beta_3, \beta_4, \beta_{5}\\
A_{4}=(\alpha_{4},\alpha_5) & \ni \beta_6, \beta_7\\
A_{5}=(\alpha_{5},\infty) & \\
&
\end{array}
\]

\noindent Therefore
\[
\GEC(F,G)=(0,1,3,2,0).
\]
Note that we counted eigenvalues of $G$ with multiplicity.

\end{example}

Next, the challenge is to extend this definition to pairs of matrices which are not generic; that is, matrices which share eigenvalues.
In the current definition, we count eigenvalues of $G$ which lie between the \textit{ open  }intervals spanned by adjacent eigenvalues of $F$; thus, if $F$ and $G$ share eigenvalues, some eigenvalues of $G$ will not be counted. Hence, we need to expand this definition for matrices which share eigenvalues.

{ 
  Consider the following example.
  \begin{example}[Eigenvalue configuration of arbitrary matrices with simple eigenvalues]
Let $F \in\mathbb{R}^{2 \times2}$ and $G \in\mathbb{R}^{3 \times3} $ be symmetric
matrices such that their corresponding eigenvalues are
\[
\alpha=(0,1),\qquad\beta=(1, 2, 3).
\]
Note that in this case $F$ and $G$ both have simple eigenvalues, but the pair is not generic since $\alpha_2 = \beta_1 = 1$.
The eigenvalues are arranged on the real line as in the following diagram.

\[
  \begin{tikzpicture}
    \draw (-1, 0) -- (4, 0); \fill [red] (0,0) circle(3pt)
    node[label=$\alpha_1$]{}; \fill [red] (1,0) circle(3pt)
    node[label=$\alpha_2$]{}; 
    \fill [blue] (1,0.8) circle(3pt)
    node[label=$\beta_1$ ]{}; \fill [blue] (2,0) circle(3pt)
    node[label=$\beta_2$ ]{}; \fill [blue] (3,0) circle(3pt)
    node[label=$\beta_3$ ]{};
  \end{tikzpicture}
  \]
  

The main difficulty here is that it is not immediately clear which interval\textemdash$(\alpha_1, \alpha_2)$ or $(\alpha_2, \infty)$\textemdash should claim $\beta_1$.
One possible solution to this problem is to ``split'' the contribution of $\beta_1$ across both of these intervals; i.e., the interval $(\alpha_1, \alpha_2)$ gets a value of $\frac12$ and the interval $(\alpha_2, \infty)$ gets a value of $\frac12$. This notion of sharing points which lie on the boundary between regions is natural and appears in many areas, such as complex analysis \cite{mariealgebraic}.

Since $\beta_2$ and $\beta_3$ both lie squarely in the interval $(\alpha_2, \infty)$, these are treated as before.
Hence, with this idea, the eigenvalue configuration is the vector
$\left( \frac{1}{2}, \frac{5}{2} \right).$

  \end{example}
  }

Let us now generalize the idea expressed in the previous example. The key observation we can make in this case is that if $(F,G)$ is not generic, then sufficiently small perturbation of their eigenvalues results in a new pair of matrices which is generic.
One natural idea { which generalizes the sharing idea from the previous example} is to define the eigenvalue configuration for arbitrary (i.e. non-generic) pairs of real symmetric matrices as the average of the eigenvalue configurations of the generic ``neighbor'' configurations, those being the configurations that can be obtained by slightly perturbing the eigenvalues.

\begin{definition}
[Eigenvalue configuration of arbitrary matrices]\label{def:EC}Let $F$~and~$G$ be real symmetric matrices. The eigenvalue configuration of~$F$ and~$G$, written as $\EC\left(  F,G\right)  $, is defined as

\begin{align*}
  \EC(F,G)  
          &= \frac{1}{2^m} \sum_{d \in \{-\varepsilon, \varepsilon\}^m} \GEC \left( F_d ,G
            \right)
\end{align*}
where 
\begin{align*}
  \varepsilon &= \text{ ``small enough'' positive number (see Remark below)} \\
  F_d &= \text{matrix whose eigenvalues $\alpha_i + d_i$ for $i = 1, \dots, m$.}
\end{align*}
\end{definition}

\begin{remark}
  \label{remark:generalEC}
  Some notes on Definition \ref{def:EC}:
  \begin{enumerate}
  \item The number $\varepsilon$ is ``small enough'' if no eigenvalue of $F$ is perturbed far enough to reach or cross the next or previous { distinct} eigenvalue { of $G$}.
    { In Section \ref{sec:arbitrary}, we give an example construction of a suitable $\varepsilon$.}
  \item { Note that any $\varepsilon$ which is small enough (in the sense above) will give the same result for $\EC(F,G)$.} 
  \item { The matrix $F_d$ is constructed by shifting the {\rm eigenvalues} of $F$, not the entries.
      Precisely, the matrix $F_d$ is constructed as
      \[
        F_d := U( \Lambda + \text{diag}(d)) U^T,
      \]
      where $F = U\Lambda U^T$ is the orthogonal eigendecomposition of (the real symmetric matrix) $F$. Note that the choice of $U$ does not matter when $F$ has repeated eigenvalues, since $\GEC$ depends only on the eigenvalues.
      We reindex the eigenvalues of $F_d$ in ascending order after perturbation.
    }
  \item Note that the pairs $(F_d, G)$ are generic for all $d \in \{-\varepsilon, \varepsilon\}^m$.
  \item If $(F, G)$ is generic, then $\GEC(F_d, G) = \GEC(F,G)$.
    This is because $\varepsilon$ is chosen so that eigenvalues of $F$ do not move far enough to cross any other eigenvalues; hence, the arrangement of the eigenvalues of generic matrices is not disturbed.
    Therefore $\EC(F,G) = \GEC(F,G)$, and so
    the generic eigenvalue configuration definition (Definition \ref{def:GEC}) is a special case of this one.
  \end{enumerate}
\end{remark}

  \begin{example}[Eigenvalue configuration of arbitrary matrices] \label{ex:g1}
    Let $F \in \mathbb{R}^{3 \times 3}$ and $G \in \mathbb{R}^{2 \times 2}$ be real symmetric matrices such that their corresponding eigenvalues are
    \[
      \alpha = (0,1,1), \qquad \beta = (1,2).
    \]
    Pictorially, the eigenvalue arrangement is as shown in
    the diagram below.
    \[
      \begin{tikzpicture}
        \draw (-1, 0) -- (3, 0);
        \fill [red] (0,0) circle(3pt) node[label=$\alpha_1$]{}; \fill [red]
        (1,0) circle(3pt) node[label=$\alpha_2$]{}; \fill [red] (1,0.75)
        circle(3pt) node[label=$\alpha_3$]{}; \fill [blue] (1,1.5)
        circle(3pt) node[label=$\beta_1$]{}; \fill [blue] (2,0)
        circle(3pt) node[label=$\beta_2$]{};
      \end{tikzpicture}
      \]
      
    Let $\varepsilon = 0.25$ (which can be verified to satisfy Definition \ref{def:EC}, as in Remark \ref{remark:generalEC}).
  Now, for each $d \in \{-\varepsilon, \varepsilon\}^3$, we compute the eigenvalue configuration of $F_d$ and $G$ in the tables below.
  Note that each pair $(F_d, G)$ is generic.
  \[
    \begin{array}{|c|c|c|}
      \hline
      \sign (d) & \text{Picture} & \GEC(F_d, G) \\
      \hline
      - - - 
                & 
                  \begin{tikzpicture}[scale=0.5]
                    \draw (-1, 0) -- (3, 0);
                    \fill [red] (-0.25,0) circle(3pt) node[]{}; \fill
                    [red] (0.75,0) circle(3pt) node[]{}; \fill [red]
                    (0.75,0.5) circle(3pt) node[]{}; \fill [blue] (1,0)
                    circle(3pt) node[]{}; \fill [blue] (2,0)
                    circle(3pt) node[]{};
                  \end{tikzpicture}
                                 &
                                   \begin{bmatrix}
                                     0 \\ 0 \\ 2
                                   \end{bmatrix} \\
      - - + 
                & 
                  \begin{tikzpicture}[scale=0.5]
                    \draw (-1, 0) -- (3, 0);
                    \fill [red] (-0.25,0) circle(3pt) node[]{}; \fill
                    [red] (0.75,0) circle(3pt) node[]{}; \fill [red]
                    (1.25,0) circle(3pt) node[]{}; \fill [blue] (1,0)
                    circle(3pt) node[]{}; \fill [blue] (2,0)
                    circle(3pt) node[]{};
                  \end{tikzpicture}
                                 &
                                   \begin{bmatrix}
                                     0 \\ 1 \\ 1
                                   \end{bmatrix} \\
      - + - 
                & 
                  \begin{tikzpicture}[scale=0.5]
                    \draw (-1, 0) -- (3, 0);
                    \fill [red] (-0.25,0) circle(3pt) node[]{}; \fill
                    [red] (1.25,0) circle(3pt) node[]{}; \fill [red]
                    (0.75,0) circle(3pt) node[]{}; \fill [blue] (1,0)
                    circle(3pt) node[]{}; \fill [blue] (2,0)
                    circle(3pt) node[]{};
                  \end{tikzpicture}
                                 &
                                   \begin{bmatrix}
                                     0 \\ 1 \\ 1
                                   \end{bmatrix} \\
      - + + 
                & 
                  \begin{tikzpicture}[scale=0.5]
                    \draw (-1, 0) -- (3, 0);
                    \fill [red] (-0.25,0) circle(3pt) node[]{}; \fill
                    [red] (1.25,0.5) circle(3pt) node[]{}; \fill [red]
                    (1.25,0) circle(3pt) node[]{}; \fill [blue] (1,0)
                    circle(3pt) node[]{}; \fill [blue] (2,0)
                    circle(3pt) node[]{};
                  \end{tikzpicture}
                                 &
                                   \begin{bmatrix}
                                     1 \\ 0 \\ 1
                                   \end{bmatrix} \\
      \hline
    \end{array} \qquad 
    \begin{array}{|c|c|c|}
      \hline
      \sign (d) & \text{Picture} & \GEC(F_d, G) \\
      \hline
      + - - 
                & 
                  \begin{tikzpicture}[scale=0.5]
                    \draw (-1, 0) -- (3, 0);
                    \fill [red] (0.25,0) circle(3pt) node[]{}; \fill
                    [red] (0.75,0) circle(3pt) node[]{}; \fill [red]
                    (0.75,0.5) circle(3pt) node[]{}; \fill [blue] (1,0)
                    circle(3pt) node[]{}; \fill [blue] (2,0)
                    circle(3pt) node[]{};
                  \end{tikzpicture}
                                 &
                                   \begin{bmatrix}
                                     0 \\ 0 \\ 2
                                   \end{bmatrix} \\
      + - + 
                & 
                  \begin{tikzpicture}[scale=0.5]
                    \draw (-1, 0) -- (3, 0);
                    \fill [red] (0.25,0) circle(3pt) node[]{}; \fill
                    [red] (0.75,0) circle(3pt) node[]{}; \fill [red]
                    (1.25,0) circle(3pt) node[]{}; \fill [blue] (1,0)
                    circle(3pt) node[]{}; \fill [blue] (2,0)
                    circle(3pt) node[]{};
                  \end{tikzpicture}
                                 &
                                   \begin{bmatrix}
                                     0 \\ 1 \\ 1
                                   \end{bmatrix} \\
      + + - 
                & 
                  \begin{tikzpicture}[scale=0.5]
                    \draw (-1, 0) -- (3, 0);
                    \fill [red] (0.25,0) circle(3pt) node[]{}; \fill
                    [red] (1.25,0) circle(3pt) node[]{}; \fill [red]
                    (0.75,0) circle(3pt) node[]{}; \fill [blue] (1,0)
                    circle(3pt) node[]{}; \fill [blue] (2,0)
                    circle(3pt) node[]{};
                  \end{tikzpicture}
                                 &
                                   \begin{bmatrix}
                                     0 \\ 1 \\ 1
                                   \end{bmatrix} \\
      + + + 
                & 
                  \begin{tikzpicture}[scale=0.5]
                    \draw (-1, 0) -- (3, 0);
                    \fill [red] (0.25,0) circle(3pt) node[]{}; \fill
                    [red] (1.25,0.5) circle(3pt) node[]{}; \fill [red]
                    (1.25,0) circle(3pt) node[]{}; \fill [blue] (1,0)
                    circle(3pt) node[]{}; \fill [blue] (2,0)
                    circle(3pt) node[]{};
                  \end{tikzpicture}
                                 &
                                   \begin{bmatrix}
                                     1 \\ 0 \\ 1
                                   \end{bmatrix} \\ \hline
    \end{array}
    \]
  Then, applying Definition \ref{def:EC}, we have
     \[ 
     \EC(F,G) \;\; =\;\; \frac{1}{2^m} \sum_{d \in \{-\varepsilon, \varepsilon\}^m} \GEC(F_d, G) 
            \;\; =\;\; \frac{1}{2^3} \left(
               2
               \begin{bmatrix}
                 0 \\ 0 \\ 2
               \end{bmatrix}
               +
               4
               \begin{bmatrix}
                 0 \\ 1 \\ 1
               \end{bmatrix}
               +
               2
               \begin{bmatrix}
                 1 \\ 0 \\ 1
               \end{bmatrix}
               \right) 
            \;\;=\;\; \frac{1}{8}
               \begin{bmatrix}
                 2 \\ 4 \\ 10 
               \end{bmatrix} 
            \;\;=\;\;
               \begin{bmatrix}
                 1/4 \\ 1/2 \\ 5/4.
               \end{bmatrix}.
   \]
\end{example}

}

\subsection{Stating the problem}

\label{sec:problem}In this section, we will state the problem precisely.
The goal of this paper is to develop an algorithm for the following problem.

\begin{problem}
\ 
\begin{enumerate}
\item[In:\ ] $F \in \mathbb{R}[p]^{m \times m}$ and  $G \in \mathbb{R}[p]^{n \times n}$,  symmetric matrices
where $p$ is a finite set of parameters.

$c \hspace{0.35em}\in\mathbb{Q}^{m}$, an eigenvalue configuration.

\item[Out:] a \textquotedblleft simple\ condition\textquotedblright\ on 
$p$ such that $c=\EC\left(F,G\right)  $.
\end{enumerate}

\end{problem}

\noindent Let us now clarify precisely what the problem is and what the
challenge is. This amounts to clarifying what
``simple\ condition'' means.
To illustrate, consider the following example. Suppose that $F$~and~$G$ are generic.
We will write down an equivalent condition for $c = \GEC(F,G)$ by explicitly writing all quantifiers and all implicit assumptions.
We have
\begin{equation}
c=\GEC\left(  F,G\right)  \ \ \ :\Longleftrightarrow\ \ \ \underset{\beta
_{1}\leq\cdots\leq\beta_{n}}{\underset{\alpha_{1}\leq\cdots\leq\alpha
_{m}}{\underset{\beta_{1},\ldots,\beta_{n}}{\underset{\alpha_{1},\ldots
,\alpha_{m}}{\exists}}}}\bigwedge\left(
\begin{array}
[c]{l}%
\underset{1\leq i\leq m}{\bigwedge}\left\vert \alpha_{i}I_{m}-F\right\vert
=0\\
\underset{x}{\forall}\left\vert xI_{m}-F\right\vert =0\Longrightarrow\left(
x=\alpha_{1}\vee\cdots\vee x=\alpha_{m}\right)  \\
  \hdashline 
\underset{1\leq i\leq n}{\bigwedge}\left\vert \beta_{i}I_{n}-{G}\right\vert =0\\
\underset{x}{\forall}\left\vert xI_{n}-{G}\right\vert \ =0\Longrightarrow\left(
x=\beta_{1}\vee\cdots\vee x=\beta_{n}\right)  \\
  \hdashline
\underset{1\leq j\leq c_{0}}{\bigwedge}\beta_{j}<\alpha_{1}\\
\underset{c_{0}+1\leq j\leq c_{0}+c_{1}}{\bigwedge}\alpha_{1}<\beta_{j}%
<\alpha_{2}\\
\hspace{1.3em}\vdots\\
\underset{c_{0}+\cdots+c_{m-1}+1\leq j\leq c_{0}+\cdots+c_{m-1}+c_{m}%
}{\bigwedge}\alpha_{m}<\beta_{j}%
\end{array}
\right)  \label{eq:qe}%
\end{equation}
where $c_{0}$ is a short hand for $n-c_{1}-\cdots-c_{m}$. The top block states
that all the eigenvalues of $F$ are exactly~$\alpha_{1},\ldots,\alpha_{m}$.
The middle block states that all the eigenvalues of $G$ are exactly $\beta
_{1},\ldots,\beta_{n}$. The bottom block states that the eigenvalue configuration
is given by $c$, and that $F$~and~$G$ are generic.

Obviously, the above is a condition on the parameters of $F$~and~$G$ so that their eigenvalues are arranged in the given way.
However, it is not ``simple'' in the sense that the expression involves quantified variables $\alpha$,~$\beta$, and $x$.

Note that the whole expression (\ref{eq:qe}) is a quantified Boolean combination of
polynomial equations and inequalities; in other words, it is a well-formed formula in
the first-order theory of real closed fields.
Alfred Tarski's celebrated theorem on real closed fields states that every
first-order formula in the theory of real closed fields is equivalent to a quantifier-\emph{free}
formula~\cite{Tarski:30}. Thus, by \textquotedblleft simple\ condition'', we
mean a quantifier-\emph{free} formula in the parameters of $F$~and~$G$;
concretely, this means an expression which can be readily rewritten as a Boolean combination of equalities and inequalities of polynomials in the parameters of $F$~and~$G$. 


{ 
  Another important consideration in what constitutes a simple condition is the size and complexity of the output. In particular, Tarski's theorem guarantees that any quantified formula over the real numbers can be written as a  quantifier-free formula in the language of the reals; precisely, this means a logical condition involving only real variables, polynomial ring operations, equality and inequality symbols, and logical conjunction and disjunction. However, such expressions can be unwieldy and difficult to understand, particularly when produced by general quantifier elimination algorithms. To mitigate this, we can introduce certain higher-level language constructs to improve the readability and reduce the size of the output condition. In this paper, we will provide an output condition which allows the higher-level matrix operations of multiplication, inverse, and signature. These reduce the size and complexity of the output significantly while still being quantifier-free.

Thus, the challenge is to produce a quantifier-free formula for
all $m$ and $n$ efficiently and in a structured form using only 
real numbers,
equality, 
matrix arithmetic,
and matrix signature.
}
\medskip

 As mentioned in the introduction, 
the proposed challenge generalizes Descartes' rule of signs, a classical real root counting method with many applications.
Hence we expect that a generalization may have even more applications.

One potential application is the reachability problem for rank-$r$ updates.
Let $F\in\mathbb{R}^{m\times m}$ be symmetric and let $r\leq m$. A positive semidefinite rank-$r$
update of $F$ is
\[
G=F+E,
\]
where $E$ is a semidefinite real symmetric matrix of rank at most $r$.
A natural question is how the eigenvalues of $F$ are affected under a rank update;  
that is, where are the relative locations of the eigenvalues of the updated matrix $G$?

{ For arbitrary $r$, Weyl's inequalities (\cite{buch2006eigenvalues}) state that $\alpha_j \le \beta_{j} \le \alpha_{j+r}$ for all $j$; however, both these conditions are conditions on the eigenvalues, and therefore involve quantifiers in the case of parametric $F$ and $G$. Our main theorem (Theorem \ref{thm:main}) provides polynomial conditions in the parameters which explicitly delineate when the eigenvalue configuration changes, and therefore provide a starting point for algebraically studying rank updates.} 

{ Another natural followup question is: for a given parametric $F$ and $G$, does there exist positive semidefinite $E$ of rank at most $r$ so that $G = F + E$?  In Section \ref{sec:update}, we give an explicit result in this direction.}
     

\section{Main Result}

\label{sec:main}

In this section, we will state our main theorem. For this, we introduce the two central notions used in the main result. The first is a \textit{combinatorial} object which depends only on the size of $F$, and the second is an \textit{algebraic} object constructed from the parameters of $F$ and $G$.

\begin{definition}[Combinatorial part]
  \label{def:C}
  The matrix $\Cs \in \mathbb{Q}^{m \times 2^m}$ is defined as
  \[
    \Cs = V_mH_m^{-1}
  \]
  where
  \begin{enumerate}
    \item
      $H_m$ is the $2^m \times 2^m$ Hadamard matrix 
      whose rows are indexed by $e \in \{0,1\}^m$ and whose columns are indexed by $s \in \{-, +\}^m$ where
      $$ (H_m)_{e,s} = s_1^{e_1} \cdots s_m ^{e_m},$$
      where $-$ is treated as $-1$ and $+$ treated as $1$ in the calculation.

      \medskip

      We will drop the subscript and call this matrix $H$ when the context is clear.
      { This is the same standard sign-determination matrix (\cite[Lecture~VII]{yap}) used in \cite{BKR1986}}.

    \item $V_m \in\mathbb{Z}^{m \times2^{m}}$ is the matrix whose rows are indexed
      by $t\in\left[  1,\ldots,m\right]  $ and whose columns are lexicographically
      indexed by $s\in\left\{  -,+\right\}  ^{m}$ where
      \[
        (V_m)_{t,s}=\one_{v(s,+),m-t}%
      \]
      where $\one$ denotes the Kronecker delta function, and $v(s,+)$ is the sign variation count of the sequence obtained by appending a $+$ to $s$; that is, the number of times two consecutive elements of $(s, +)$ are opposite sign, ignoring zeros.

      \medskip

      We will drop the subscript and call this matrix $V$ when the context is clear.
  \end{enumerate}
\end{definition}

{ 
  \begin{remark}
    Note that since $H_m$ is a Hadamard matrix, we have that
    \begin{align*}
      H_m^{-1} = \frac{1}{2^m} H_m^T.
    \end{align*}
    Hence, in practice, it is not necessary to compute any inverses in the computation of $\Cs$.
    In this paper, we make the decision to leave $H_m^{-1}$ in the above definition to emphasize its conceptual importance.
  \end{remark}
}

\begin{example} \label{ex:C}
  Let $m = 2$. We will produce the matrix $\Cs$.
  \begin{enumerate}
    \item First, we construct $H_2$.
      \[
        H_2=
        \left[
          \begin{array}{c|cccc}
            e \backslash s & - - & -+ & +- & ++ \\
            \hline
            00 & 1 & 1 & 1 & 1 \\
            01 & -1 & 1 & -1 & 1 \\
            10 & -1 & -1 & 1 & 1 \\
            11 & 1 & -1 & -1 & 1 \\
          \end{array} \right]
      \]

      Explanations on a few entries:%
      \begin{align*}
        (H_2)_{01, - -}  &  =(-1)^0(-1)^1 = -1 \\
        (H_2)_{00, - +}  &  =(-1)^0(1)^0 = 1 \\
        (H_2)_{11, + -}  &  =(1)^1(-1)^1 = -1. 
      \end{align*}

    \item Then, we construct $V_2.$
      \[
        V_2=
        \left[
          \begin{array}
            [c]{l|cccc}%
            t\backslash s & -- & -+ & +- & ++\\\hline
            1 & 1 & 1 & 0 & 0\\
            2 & 0 & 0 & 0 & 1
          \end{array}
        \right]
      \]

      Explanations on a few entries:%
      \begin{align*}
        V_{1,--}  &  =\one_{v(--+),2-1}=\one_{1,1}=1\\
        V_{2,--}  &  =\one_{v(--+),2-2}=\one_{1,0}=0\\
        V_{2,++}  &  =\one_{v(+++),2-2}=\one_{0,0}=1 .
      \end{align*}
  \end{enumerate}
  Putting these together, we get that
  \begin{align*}
    \Cs &= V_2H_2^{-1}  
            = 
        \left[
          \begin{array}
            [c]{cccc}%
             1 & 1 & 0 & 0\\
             0 & 0 & 0 & 1
          \end{array}
              \right]
        \left[
          \begin{array}{cccc}
             1 & 1 & 1 & 1 \\
             -1 & 1 & -1 & 1 \\
             -1 & -1 & 1 & 1 \\
             1 & -1 & -1 & 1 \\
          \end{array} \right]^{-1} 
    = \left[\begin{array}{cccc}
\frac{1}{2} & 0 & -\frac{1}{2} & 0 
\\
 \frac{1}{4} & \frac{1}{4} & \frac{1}{4} & \frac{1}{4} 
\end{array}\right]. 
  \end{align*} 
\end{example}

\begin{definition}[Algebraic part] \label{def:A}
  $\As$, or $\As(F,G)$ to emphasize the dependence on the parameters of $F$ and $G$, is the column vector of length $2^m$ whose entries are indexed by $e \in \{0, 1\}^m$ where
  \[
    \left(\As\right)_e = \sigma( f_e(G) )
  \]
  where
  \begin{itemize}
  \item 
  $\sigma$ denotes the signature of a matrix (i.e., the difference between the number of positive eigenvalues and number of negative eigenvalues)
\item 
      $f$ is the characteristic polynomial of $F$
    \item $f_e$ is the polynomial
      \begin{align*}
        f_e &= \left( f^{(0)} \right)^{e_0} \cdots \left( f^{(m-1)} \right)^{e_{m-1}},
      \end{align*}
      with $f^{(k)}$ denoting the $k$-th derivative of $f$. 
  \end{itemize}
\end{definition}

\begin{example} \label{ex:A}
    Let $m=2$ and $n=3$. Then let $F$ and $G$ be matrices where each entry is an independent parameter; that is,
    \[
      F=%
      \begin{bmatrix}
        a_{1,1} & a_{1,2}\\
        a_{1,2} & a_{2,2}%
      \end{bmatrix}
      ,\qquad G=%
      \begin{bmatrix}
        b_{1,1} & b_{1,2} & b_{1,3} \\
        b_{1,2} & b_{2,2} & b_{2,3} \\
        b_{1,3} & b_{2,3} & b_{3,3}
      \end{bmatrix}
      .
    \]
    We will compute one of the entries of $\As(F,G)$. By Definition \ref{def:A} we have that
    \[
      \As =
      \begin{bmatrix}
        \sigma( f_{00}(G) ) \\
        \sigma( f_{01}(G) ) \\
        \sigma( f_{10}(G) ) \\
        \sigma( f_{11}(G) )
      \end{bmatrix}.
    \]
    As an example, consider $e = (0,1)$. We will compute $f_{01}(G)$.
    We have
    \begin{align*}
      f  &  =\det(xI_{2}-F)=x^{2}+(-a_{2,2}-a_{1,1})x+a_{2,2}a_{1,1}-a_{1,2}^{2}\\
      f_{(0,1)}  &  =f^{(0)^{0}}f^{(1)^{1}}
           =2x-a_{2,2}-a_{1,1}.
    \end{align*}
    Next, we compute
    \begin{align*}
      f_{(0,1)}(G)  &  =2G+(-a_{2,2}-a_{1,1})I_{3}\\
                    &  =%
      \left[ \begin{array}{ccc}
        -a_{2,2}-a_{1,1}+2 b_{1,1} & 2 b_{1,2} & 2 b_{1,3} 
        \\
        2 b_{1,2} & -a_{2,2}-a_{1,1}+2 b_{2,2} & 2 b_{2,3} 
        \\
        2 b_{1,3} & 2 b_{2,3} & -a_{2,2}-a_{1,1}+2 b_{3,3} 
      \end{array}\right ]
    \end{align*}
    Finally, we have that
    \begin{align*}
      \sigma(f_{01}(G)) = \# \text{positive eigenvalues of $f_{01}(G)$} \,\, - \,\, \# \text{negative eigenvalues of $f_{01}(G)$},
    \end{align*}
    which clearly depends on the particular values of the parameters of $F$ and $G$. 
\end{example}

\bigskip

\noindent Now, we are ready to state the main theorem which relates these two notions.
{ 
  \begin{theorem}
    [Main Result]\label{thm:main}Let $F\in\mathbb{R}^{m\times m}$ and
    $G\in\mathbb{R}^{n\times n}$ be real symmetric matrices. We have%
    \[
      \EC\left(  F,G\right) \,\,\;  = \; \,\,  \Cs\;\;\As(F,G).
    \]

  \end{theorem}

  \begin{remark}
    Note that the entries of the symmetric matrices $f_e(G)$, whose signatures make up the vector $\As(F,G)$, are polynomials in the entries $a_{ij}$ of~$F$ and $b_{ij}$ of $G$. 
    Hence, their signatures can be written as {quantifier-free} conditions on polynomials in the parameters of~$F$~and~$G$ since Descartes' rule of signs on the characteristic polynomial gives exact counts of the numbers of positive and negative eigenvalues. Thus the condition is quantifier-free. (Expanding the signatures into sign conditions yields nontrivial output size, and we leave the study of their structure for a future work.)
  \end{remark}
  \begin{remark}
    In case the reader is familiar with our related work \cite{symmetric}, the reader will notice that the main theorem in that paper is stated almost identically to Theorem \ref{thm:main} in this paper. However, they are completely different results, as the respective $C$ (combinatorial parts) and $A$ (algebraic parts) are defined as completely different concepts in each paper. The main theorem in this paper is based on the signature of matrices while the main theorem in \cite{symmetric} is based on symmetric polynomials.
    However, in both theorems, there is a combinatorial part and an algebraic part; hence, we structure the theorem statements similarly to emphasize the super-structural similarities. 
  \end{remark}

  \begin{example}
    [General case with parametric matrices]
    \label{ex:general}
    Let $m = 2$ and $n = 3$ and let $F$ and $G$ be as in Example \ref{ex:A}; that is, 
    \[
      F =
      \begin{bmatrix}
        a_{11} & a_{12} \\
        a_{12} & a_{22} \\
      \end{bmatrix} , \qquad \qquad
      G =
      \begin{bmatrix}
        b_{11} & b_{12} & b_{13} \\
        b_{12} & b_{22} & b_{23} \\
        b_{13} & b_{23} & b_{33}
      \end{bmatrix},
    \]
    where each $a_{ij}$ and $b_{ij}$ is a fully independent parameter.
    We will now use Theorem \ref{thm:main} to write a condition on the parameters of $F$ and $G$ which is equivalent to the eigenvalues of $F$ and $G$ being arranged as in
    the following diagram.
    \[
      \begin{tikzpicture}
        \draw (-1, 0) -- (6, 0); \fill [red] (3,0) circle(3pt)
        node[label=$\alpha_1$]{}; \fill [red] (5,0) circle(3pt)
        node[label=$\alpha_2$]{}; \fill [blue] (4,0) circle(3pt)
        node[label=$\beta_3$]{}; \fill [blue] (1,0) circle(3pt)
        node[label=$\beta_2$]{}; \fill [blue] (0,0) circle(3pt)
        node[label=$\beta_1$]{};
      \end{tikzpicture}
      \]
      
    That is, we will find a condition for $\EC(F,G) =
    \begin{bmatrix}
      1 \\ 0
    \end{bmatrix}.$
    By Theorem \ref{thm:main}, we have that
    \begin{align*}
      \EC(F,G) =
      \begin{bmatrix}
        1 \\ 0
      \end{bmatrix}
      \,\, \iff \,\,
      \begin{bmatrix}
        1 \\ 0
      \end{bmatrix}
      &= \Cs \,\, \As(F,G). 
    \end{align*}
    In Example \ref{ex:C}, we found that
    \[
      \Cs = 
    \left[\begin{array}{cccc}
\frac{1}{2} & 0 & -\frac{1}{2} & 0 
\\
 \frac{1}{4} & \frac{1}{4} & \frac{1}{4} & \frac{1}{4} 
\end{array}\right],
    \]
    and from Example \ref{ex:A} we have
    \[
      \As(F,G) = 
      \begin{bmatrix}
        \sigma( f_{00}(G) ) \\
        \sigma( f_{01}(G) ) \\
        \sigma( f_{10}(G) ) \\
        \sigma( f_{11}(G) )
      \end{bmatrix}.
    \]
    Thus by Theorem \ref{thm:main} we have
    \begin{align*}
      \EC(F,G) =
      \begin{bmatrix}
        1 \\ 0
      \end{bmatrix}
      \,\, \iff \,\,
      \begin{bmatrix}
        1 \\ 0
      \end{bmatrix}
      &=  
    \left[\begin{array}{cccc}
\frac{1}{2} & 0 & -\frac{1}{2} & 0 
\\
 \frac{1}{4} & \frac{1}{4} & \frac{1}{4} & \frac{1}{4} 
\end{array}\right]
      \begin{bmatrix}
        \sigma( f_{00}(G) ) \\
        \sigma( f_{01}(G) ) \\
        \sigma( f_{10}(G) ) \\
        \sigma( f_{11}(G) )
      \end{bmatrix}.
    \end{align*}
    With the fact that $-n \le \sigma(f_e(G)) \le n$ and $\sigma(f_{00}(G)) = n$ always, one can brute force the solutions to the above linear system and find that $EC(F,G) =
    \begin{bmatrix}
      1 \\ 0
    \end{bmatrix}$ is equivalent to
    \begin{align*}
      &\begin{bmatrix}
        \sigma( f_{00}(G) ) \\
        \sigma( f_{01}(G) ) \\
        \sigma( f_{10}(G) ) \\
        \sigma( f_{11}(G) )
      \end{bmatrix}
      \in
      \left\{
      \begin{bmatrix}
        3\\-3\\1\\-1
      \end{bmatrix},
      \begin{bmatrix}
        3\\-2\\1\\-2
      \end{bmatrix},
      \begin{bmatrix}
        3\\-1\\1\\-3
      \end{bmatrix}
      \right\}.
    \end{align*}
    Note that all the work done in this example to enumerate the possible signature vectors used only $m$ and $n$ and not the entries of the matrices; thus, the above set of possible signature vectors could be pre-computed for fixed $m$ and $n$ and used to construct a  condition for {\rm any} $m \times m$ matrix $F$ and {\rm any} $n \times n$ matrix $G$ having the same eigenvalue configuration.
    
  \end{example}

  \begin{example}[Matrices with only one parameter]
    In Example \ref{ex:general}, we studied $F$ and $G$ which were symmetric matrices, each of whose elements was an independent parameter.
    However, it need not be the case that each entry is an independent parameter. We would, for instance, have the entries of both matrices be polynomials in some smaller list of parameters.
    Consider, for example, the matrices
    \begin{align*}
      F =
      \begin{bmatrix}
        1 & -2 \\ -2 & 1 
      \end{bmatrix}
    \end{align*}
    and
    \begin{align*}
      G =
      \left[\begin{array}{ccc}
        1 & p^{2}+p  & p +1 
        \\
        p^{2}+p  & -p  & p -1 
        \\
        p +1 & p -1 & p  
      \end{array}\right],
    \end{align*}
    where $p$ is a (single) parameter.
    Suppose we wish to know for which values of $p$ the eigenvalues of the matrix $G$ interlace with that of $F$; that is, when $\EC(F, G) = (1,1)$.
    By Theorem \ref{thm:main} and reusing the computation from Example \ref{ex:C}, we have that
    \begin{align}
      \EC(F,G) =
      \begin{bmatrix}
        1 \\ 1
      \end{bmatrix} \, \, &\iff \,\,
                            \begin{bmatrix}
                              1 \\ 1 
                            \end{bmatrix} =
                            \Cs \,\, \As(F,G) \nonumber \\
    &\iff \,\, 
    \begin{bmatrix}
      1 \\ 1 
    \end{bmatrix} =
\left[\begin{array}{cccc}
\frac{1}{2} & 0 & -\frac{1}{2} & 0 
\\
 \frac{1}{4} & \frac{1}{4} & \frac{1}{4} & \frac{1}{4} 
\end{array}\right]
\left[
\begin{array}
[c]{c}%
\sigma(f_{00}(G)) \\
\sigma(f_{01}(G)) \\
\sigma(f_{10}(G)) \\
\sigma(f_{11}(G)) \\
\end{array}
\right], \label{eq:example}
  \end{align}
  where
  \begin{align*}
    \sigma(f_{00}(G)) &= \sigma \left( \left[\begin{array}{ccc}
1 & 0 & 0 
\\
 0 & 1 & 0 
\\
 0 & 0 & 1 
    \end{array}\right] \right)\\
    \sigma(f_{01}(G)) &= \sigma\left(\left[\begin{array}{ccc}
0 & 2 p^{2}+2 p  & 2 p +2 
\\
 2 p^{2}+2 p  & -2-2 p  & 2 p -2 
\\
 2 p +2 & 2 p -2 & 2 p -2 
    \end{array}\right] \right) \\
    \sigma(f_{10}(G)) &= \sigma\left(\left[\begin{array}{ccc}
p^{4}+2 p^{3}+2 p^{2}+2 p -3 & -p^{3}-p^{2}-p -1 & \left(p -1\right) \left(p +1\right)^{2} 
\\
 -p^{3}-p^{2}-p -1 & p^{4}+2 p^{3}+3 p^{2}-2 & p^{3}+2 p^{2}-p +2 
\\
 \left(p -1\right) \left(p +1\right)^{2} & p^{3}+2 p^{2}-p +2 & 3 p^{2}-2 p -1 
    \end{array}\right] \right)\\
    \sigma(f_{11}(G)) &=  \sigma\left(\scalebox{0.75}{$\left[\begin{array}{ccc}
-2 p^{5}-2 p^{4}-4 p^{2}-6 p -2 & 2 p^{6}+6 p^{5}+12 p^{4}+12 p^{3}-2 p^{2}-2 p +4 & 2 \left(p^{3}+p^{2}+p -1\right) \left(p +1\right)^{2} 
\\
 2 p^{6}+6 p^{5}+12 p^{4}+12 p^{3}-2 p^{2}-2 p +4 & -4 p^{5}-8 p^{4}-12 p^{3}-16 p^{2}+8 p  & 2 p^{5}+2 p^{4}-16 p^{2}-2 p -2 
\\
 2 \left(p^{3}+p^{2}+p -1\right) \left(p +1\right)^{2} & 2 p^{5}+2 p^{4}-16 p^{2}-2 p -2 & 4 \left(p -1\right) \left(p^{3}+4 p^{2}+1\right) 
    \end{array}\right]$}\right).
  \end{align*}
      Taking into account the fact that $G$ is a $3\times3$ matrix, and so $-3 \le \sigma(f_e(G)) \le 3$, and the fact that $\sigma(f_{00}(G)) = n = 3$, one can find the (finitely many) integer vector solutions to the linear system by brute force and find that (\ref{eq:example}) is equivalent to 
    \[
      \begin{bmatrix}
        \sigma(f_{00}(G)) \\ \sigma(f_{01}(G)) \\ \sigma(f_{10}(G)) \\ \sigma(f_{11}(G))
      \end{bmatrix}
      \in
      \left\{
        \begin{bmatrix}
          3\\-3\\1\\3
        \end{bmatrix}, 
        \begin{bmatrix}
          3\\-2\\1\\2
        \end{bmatrix}, 
        \begin{bmatrix}
          3\\-1\\1\\1
        \end{bmatrix}, 
        \begin{bmatrix}
          3\\0\\1\\0
        \end{bmatrix},
        \begin{bmatrix}
          3\\1\\1\\-1
        \end{bmatrix},
        \begin{bmatrix}
          3\\2\\1\\-2
        \end{bmatrix},
        \begin{bmatrix}
          3\\3\\1\\-3
        \end{bmatrix}
      \right\}.
    \]

Note that the above is a condition which depends only on the parameter $p$ and not on the eigenvalues of $F$~and~$G$.
  As stated previously, the signatures of each (real symmetric) matrix $f_e(G)$ can be expanded using Descartes' rule of signs to express the positive and negative eigenvalues of $f_e(G)$ in terms of the sign variation count of the respective characteristic polynomials.

  For example, consider
  $$f_{01}(G) = \left[\begin{array}{ccc}
0 & 2 p^{2}+2 p  & 2 p +2 
\\
 2 p^{2}+2 p  & -2-2 p  & 2 p -2 
\\
 2 p +2 & 2 p -2 & 2 p -2 
    \end{array}\right].$$
  We have
  \begin{align*}
    \operatorname*{charpoly} f_{01}(G) = x^3 + 4x^2 - (4p^4 + 8p^3 + 16p^2 + 4) x + 8p^5 - 8p^4 - 32p^3 - 16p^2 - 8p - 8.
  \end{align*}
  Since $f_e(G)$ is real symmetric, its characteristic polynomial has all real roots; thus by Descartes' rule of signs, the number of positive eigenvalues of $f_{01}(G)$ equals
  \begin{align*}
    v(1, 4, -(4p^4 + 8p^3 + 16p^2 + 4), 8p^5 - 8p^4 - 32p^3 - 16p^2 - 8p - 8)
  \end{align*}
  and the number of negative eigenvalues equals
  \begin{align*}
    v(-1, 4, (4p^4 + 8p^3 + 16p^2 + 4), 8p^5 - 8p^4 - 32p^3 - 16p^2 - 8p - 8).
  \end{align*}
  Therefore
  \begin{align*}
    \sigma(f_{01}(G)) = 
    \,\,\,\,\,\,&v(1, 4, -(4p^4 + 8p^3 + 16p^2 + 4), 8p^5 - 8p^4 - 32p^3 - 16p^2 - 8p - 8) \\
    - &v(-1, 4, (4p^4 + 8p^3 + 16p^2 + 4), 8p^5 - 8p^4 - 32p^3 - 16p^2 - 8p - 8).
  \end{align*}
  This expression can then be expanded to explicitly enumerate the possible signs for each of the polynomials, and the process can be repeated for $\sigma(f_{00}(G)), \sigma(f_{10}(G)),$ and $\sigma(f_{11}(G))$ to obtain the full condition involving only $p$.
  
\end{example}


  
}
{In pseudocode form, Theorem \ref{thm:main} can be restated as the following algorithm.

\begin{algorithm}
[Main Theorem: Condition for EC]\ \label{alg:sig}

\begin{itemize}
\item[In: ] $F \in \mathbb{R}[p]^{m \times m}$ and  $G \in \mathbb{R}[p]^{n \times n}$, symmetric, where $p$ is a finite list of parameters

  $c \in \mathbb{Q}^m$, eigenvalue configuration

\item[Out:] $P,$ quantifier-free formula on $p$ such that $c=\EC\left(  F,G\right)$ 
\end{itemize}

\begin{enumerate}
\item $f \longleftarrow \det(xI_m - F)$
\item For $e \in \{0,1\}^m$ do

  \medskip

   \qquad $f_e \leftarrow \left( f^{(0)}\right)^{e_0} \cdots \left( f^{(m-1)}\right)^{e_{m-1}}$

\item $\As \leftarrow $ vector such that $(\As)_e = \sigma(f_e(G))$ for $e \in \{0,1\}^m$

\item $V \leftarrow $ matrix such that $V_{t,s} = \delta_{v(s,+) = m-t}$ for $t \in [m]$ and $s \in \{-,+\}^m$
\item $H \leftarrow $ matrix such that $H_{e,s} = s_1^{e_1} \cdots s_m^{e_m}$ for $e \in \{0,1\}^m$ and $s \in \{-,+\}^m$
\item $\Cs \leftarrow VH^{-1}$
\item $\displaystyle P\longleftarrow \bigvee_{\substack{S \in \{-n, \dots, n\}^{2^m} \\ S\,\,:\,\, c \,=\, \Cs\,\, S}} \As = S$ \label{alg:conj}
\item Return $P$
\end{enumerate}
\end{algorithm}
}

{\begin{remark}\ 
\begin{enumerate}
\item Note that we allow the signature symbol $\sigma$ in the output formula $P$ in order to write a more compact and readable condition, as typically done in the literature (e.g. \cite{GV1999}).
\item Note that $\Cs$ only depends on the size of $F$, 
not on the size of $G$, and not on the entries of $F$ and $G$.
Thus when one works with many $(F,G)$ with same size of $F$, then $\Cs$ can be pre-computed once. 
\item We can show that the algorithm ultimately involves  $n 2^m$ polynomials in the parameters $p$. 
Hence, in practice, one should assign the smaller matrix as $F$ and the larger as $G$ (keeping in mind that $\EC(F,G) \ne \EC(G,F)$ in general).
\item Fully enumerating all terms in the conjunction in line \ref{alg:conj} still involves $O((2n+1)^{2^m})$ time; however, due to the fact that $\sigma(f_{00}(G)) = n$ trivially, most of the enumerated $S$ vectors are infeasible. We leave the pruning of the infeasible signature vectors for a future work.
\end{enumerate}
\end{remark}
}

\section{Proof / Derivation}

\label{sec:proof}

In this section, we will provide a proof for the main result (Theorem
\ref{thm:main}). We will do this in a series of lemmas that explain the motivation and reasoning behind each step. 
{ 
We divide the proof into three stages, where each stage builds from the previous.
  \begin{enumerate}
  \item In Section \ref{sec:strong}, we prove a special case of Theorem \ref{thm:main} in the case where $F$~and~$G$ are ``strongly generic;'' that is, the pair satisfies an additional condition (see below) in addition to being generic.
  \item In Section \ref{sec:generic}, we prove a slightly more general case of Theorem \ref{thm:main} in the case where $F$~and~$G$ are generic.
  \item Finally, in Section \ref{sec:arbitrary}, we prove Theorem \ref{thm:main} in full generality.
  \end{enumerate}
}


\begin{definition}
[Strongly generic]\label{def:strong} We say that the pair $(F,G)$ of real symmetric matrices is
\textbf{strongly generic} if both of the following hold.

\begin{itemize}
\item The pair $(F,G)$ is generic, i.e. $F$ and $G$ do not share eigenvalues.

\item Let $f(x)$ be the characteristic polynomial of $F$.
  For all eigenvalues $\beta$ of $G$ and for all $k \in \{ 1, \dots, m-1\}$, we have
  $ f^{(k)}(\beta) \ne 0,$
  where $f^{(k)}$ denotes the $k$-th derivative of $f$.
\end{itemize}
\end{definition}

\noindent We will illustrate the above definition
using two examples: one for strongly generic and the other for  generic but not strongly generic.
 These will be used as running examples in the following subsections.  
\begin{example}[Strongly generic] \label{ex:strong} We will be using this example as a running example throughout 
Section~\ref{sec:strong}. \medskip

\noindent Let
$F =
\begin{bmatrix}
4 & 0 \\
0 & 4 
\end{bmatrix}
\in\mathbb{R}^{2 \times2}$
and
$ G =
\begin{bmatrix}
2 & 0 & 0\\
0 & 2 & 0\\
0 & 0 & 8
\end{bmatrix}
\in\mathbb{R}^{3 \times3}.$
Then their respective eigenvalues are
$\alpha= (4, 4)$ and $\beta= (2,2,8)$
as shown in the following diagram.
\[
  \begin{tikzpicture}
    \draw (0, 0) -- (9, 0); \fill [red] (4,0) circle(3pt)
    node[label=$\alpha_1$]{}; \fill [red] (4,0.8) circle(3pt)
    node[label=$\alpha_2$]{}; \fill [blue] (2,0) circle(3pt)
    node[label=$\beta_1$]{}; \fill [blue] (2,0.8) circle(3pt)
    node[label=$\beta_2$]{}; \fill [blue] (8,0) circle(3pt)
    node[label=$\beta_3$]{};
  \end{tikzpicture}
  \]
  

\noindent So we have
$\GEC(F,G) = (c_{1}, c_{2}) = (0, 1).$
Further, we have $f = \det(xI - F) = (x-4)^2$ and $f^{(1)} = 2(x-4)$.
Note that $F$~and~$G$ are generic because they do not share any eigenvalues.
Furthermore, the derivative $f^{(1)}$ is nonzero on all eigenvalues of $G$.
Thus $F$ and
$G$ are strongly generic.

\end{example}

\begin{example}[Not strongly generic]\label{ex:generic} We will be using this example 
as a running example throughout Section~\ref{sec:generic}.

\medskip\noindent Let
$F=%
\begin{bmatrix}
1 & 0 \\
0 & 5 \\
\end{bmatrix}
\in\mathbb{R}^{2\times2}$ and $G=
\begin{bmatrix}
3 & 0 & 0\\
0 & 3 & 0\\
0 & 0 & 6
\end{bmatrix}
\in\mathbb{R}^{3\times3},$
and so their eigenvalues are
$\alpha=(1,5)$ $\beta=(3,3,6).$
Hence their eigenvalue configuration is shown in
the diagram below.

\[
  \begin{tikzpicture}
    \draw (0, 0) -- (7, 0); \fill [red] (1,0) circle(3pt)
    node[label=$\alpha_1$]{}; \fill [red] (5,0) circle(3pt)
    node[label=$\alpha_2$]{}; \fill [blue] (3,0) circle(3pt)
    node[label=$\beta_1$]{}; \fill [blue] (3,0.8) circle(3pt)
    node[label=$\beta_2$]{}; \fill [blue] (6,0) circle(3pt)
    node[label=$\beta_3$]{};
  \end{tikzpicture}
  \]
  
\noindent and so%
$\GEC(F,G)=(2,1).$
Note that $F$~and~$G$ are generic because they do not share any eigenvalues. However, note that
$f    =\det(xI-F)=(x-1)(x-5)$ and 
$f^{(1)}    = 2(x-3). $
Then we have
$f^{(1)}(\beta_{1}) = f^{(1)}(\beta_{2}) =0,$
and so $F$~and~$G$ are not strongly generic.

\end{example}


\subsection{Proof for strongly generic pair of matrices}
\label{sec:strong}

In this subsection, we first prove a special case of the main result (Theorem \ref{thm:main}) for strongly generic $F$~and~$G$. 
\begin{lemma}
\label{lem:strong} Let $F\in\mathbb{R}^{m\times m}$ and $G\in
\mathbb{R}^{n\times n}$ be a real symmetric \textbf{strongly generic} pair of matrices. Then
\[
\GEC(F,G)\,\,=\,\,\Cs\,\,\As\left(  F,G\right).
\]
\end{lemma} 
{Note that unless stated otherwise, all $F$ and $G$ in this section are assumed to be strongly generic.}
Our process of proving Lemma \ref{lem:strong} will follow the diagram below.
In the lemmas in this section, we repeatedly rewrite the statement $c = \GEC(F,G)$ with the goal of eliminating all references to $\alpha$ and $\beta$, so that we end with an expression involving only the entries of $F$~and~$G$.
  \begin{table}[H]
    \centering
    \begin{tabular}{|c|l|}
      \hline
      &\\
     $c = \GEC(F,G)$ & \\[2ex]
      \hspace{6.8em}$\Big \Vert $ by Definition \ref{def:GEC} & \\[2ex]
$c_t = \#\{j : \beta_j \in A_t \}$
      & (Expression involving $\alpha$ and $\beta$) \\[2ex]
      \hspace{4.8em}$\Big \Vert$ Lemma \ref{lem:elim_alpha} & \\[2ex]
    $c_t = \# \{j : \overline\# \{x : f(x) = 0 \wedge x > \beta_j\} = m-t\}$ & (Expression involving $\beta$) \\[2ex]
      \hspace{8.4em}$\Big \Vert$ Lemmas \ref{lem:elim_inner_count}, \ref{lem:elim_sv}, \ref{lem:elim_count_N} & \\[2ex]
    $c = \Cs\,\,\As(F,G)$ & (Expression involving only $a_{ij}$ and $b_{ij}$ and signature) \\[2ex]
      \hline
    \end{tabular}
  \end{table}

  \noindent As the first step, we will now find an expression for $\GEC(F,G)$ which does not refer to $\alpha$.

\begin{lemma}
[Eliminate $\alpha$]\label{lem:elim_alpha}
{Let $F$ and $G$ be strongly generic real symmetric matrices.} Let $c = \GEC(F,G)$. Then for all $t = 1, \dots, m$ we have

\[
c_{t}=\#\left\{  j: \overline\# \left\{  x:f\left(  x\right)
=0\ \wedge\ x>\beta_{j}\right\}  =m-t\right\}
\]
where $f\left(  x\right)  =\det(xI_{m}-F)$ is the
characteristic polynomial of $F$, and the symbol $\overline\#$ means to count with multiplicity. 
\end{lemma}

\begin{example}
[Running, Section \ref{sec:strong}]Recall that in the running example (Example
\ref{ex:strong}) we had that $c = (0,1)$.
On the other hand,
  

\[%
\begin{array}
[c]{lllllll}%
\#\{j:\overline\#\{x:f(x)=0\wedge x>\beta_{j}\}=2-1\} & = & \#\emptyset         & = & 0 & = & c_{1}\\
\#\{j:\overline\#\{x:f(x)=0\wedge x>\beta_{j}\}=2-2\} & = & \#\{3\} & = & 1 &  = & c_{2}.
\end{array}
\]

\end{example}

\begin{proof}[Proof of Lemma \ref{lem:elim_alpha}]
Recall the eigenvalue configuration condition (Definition \ref{def:GEC}):%
\begin{equation}
c_{t}=\#\{j:\beta_{j}\in A_{t}\}. \label{eq:GEC}%
\end{equation}
We will rewrite the expression on the right-hand side, with the goal of eventually reducing the problem to a set of real root counting
problems. First, we will repeatedly rewrite the condition $\beta_{j}\in
A_{t}$ with the goal of replacing the \emph{membership} checking of $\beta
_{j}$ with the \emph{counting} of~$\alpha_{i}$.
\begin{align*}
\beta_{j}\in A_{t}
 & \Longleftrightarrow \alpha_{t}<\beta_{j}<\alpha_{t+1} 
   \;\;\; \text{from the definition of} A_t \\
 & \Longleftrightarrow  \alpha_{1}\le\cdots\le\alpha_{t}<\beta_{j}<\alpha_{t+1}\le\cdots\le\alpha 
   \;\;\; \text{from the assumption that }\alpha_{1}\le\alpha_{2}\le\cdots\le\alpha_{m}\text{{}}\\
  &\Longleftrightarrow  \#\left\{  i:\alpha_{i}>\beta_{j}\right\}  =m-t 
   \;\;\; \text{from counting }i\ \text{such that }\alpha_{i}\text{\ is greater than }\beta_{j}
\end{align*}
Note that the counting of \textit{indices} of the eigenvalues $\alpha$ effectively counts the eigenvalues with multiplicity.

\medskip

\noindent By replacing
$\beta_{j}\in A_{t}$ with $\#\left\{  i:\alpha_{i}>\beta_{j}\right\}  =m-t$ in
the right-hand side of (\ref{eq:GEC}), we arrive at
\begin{align} \label{eq:elim_at}
   c_t =\#\left\{  j:\#\left\{  i:\alpha_{i}>\beta_{j}\right\} =m-t\right\}.
\end{align}


\noindent Finally, the claim of the lemma follows immediately by replacing
$\#\left\{  i:\alpha_{i}>\beta_{j}\right\}$ with $ \overline\#\left\{  x:f\left(
x\right)  =0\ \wedge\ x>\beta_{j}\right\}$ in the right-hand side of
(\ref{eq:elim_at}), where again $\overline{\#}$ means to count with multiplicity.
\end{proof}

\bigskip

At this point, we now have an expression for each component of $\GEC(F,G)$ which does not contain any~$\alpha$'s.
The next step is to eliminate the $\beta$'s.
This is somewhat more complicated than eliminating the $\alpha$'s, so we split the next stage into several steps.

\begin{enumerate}
\item In Lemma \ref{lem:elim_inner_count}, we will first eliminate the counting of roots of $f$; i.e., we find another way to express the quantity $\overline\#\{x: f(x) = 0 \wedge x > \beta_j\}$.
  This will be achieved by Descartes' rule of signs.
  With this, we will have rewritten each of the $m$ entries of $\GEC(F,G)$ as a root counting problem which counts eigenvalues of $G$ subject to a sign condition on the derivatives of $f$.
  
\item In Lemma \ref{lem:elim_sv}, we will then rewrite in matrix form. This is a necessary prerequisite for the next step, in which we solve all $m$ root counting problems together. 

\item Finally, in Lemma \ref{lem:elim_count_N}, we will solve the root counting problems by rewriting them in terms of the signature (i.e., the difference between the number of positive and negative eigenvalues) of matrices derived from $F$~and~$G$. { To do this, we will apply a root counting method from Ben-Or, Kozen, and Reif \cite{BKR1986}.}
\end{enumerate}

\bigskip

\noindent We begin by eliminating the reference to counting the roots of $f$. First, we need a notation.

\begin{notation}
The sign sequence of $f = \det(xI_m - F)$ evaluated at $x \in\mathbb{R}$, denoted $\sseq
(x)$, is
\[
\sseq(x) = \operatorname*{sign} \left(  f^{(0)} (x), \dots, f^{(m-1)}(x)
\right)  .
\]
Note that we only go to the $(m-1)$th derivative, since by construction the $m$-th derivative is always a positive constant; hence we can omit it for simpler presentation.
\end{notation}

\begin{lemma}
  [Eliminate counting roots of $f$]\label{lem:elim_inner_count}
{Let $F$ and $G$ be strongly generic real symmetric matrices.}
  Let $c = \GEC(F,G)$. For all $t = 1, \dots, m$, we have
\[
c_{t}=\#\left\{  j:v\left(  \sseq(\beta_{j}) ,+\right)
=m-t\right\}  .
\]

\end{lemma}

\begin{example}
[Running, Section \ref{sec:strong}]\label{ex:prev1} Recall that in the running
example (Example \ref{ex:strong}) we had $c=(0,1)$ and 
$\alpha= (4, 4)$ and $\beta= (2,2,8).$ We also had that
$f   =\det(xI-F)=(x-4)^2$ and 
$f^{(1)}    = 2(x-4).$

On the other hand, we can compute
\begin{align*}
v(\sseq(\beta_{1}), +) = v(\sseq(2), +)  &  = v\left(  \operatorname*{sign}\left(  f(2), f^{(1)}(2) \right)  , +\right) 
 = v(+, -, +) = 2\\
v(\sseq(\beta_{2}), +) = v(\sseq(2), +)  &  = v\left(  \operatorname*{sign}\left(  f(2), f^{(1)}(2) \right)  , +\right) = v(+, -, +) = 2\\
v(\sseq(\beta_{3}), +) = v(\sseq(8), +)  &  = v\left(  \operatorname*{sign}\left(  f(8), f^{(1)}(8) \right)  , +\right)  = v(+, +, + ) = 0.
\end{align*}
\noindent So then
\[%
\begin{array}
[c]{lllllll}%
\#\{j:v(\sseq(\beta_{j}))=2-1\} & = & \#\emptyset & = & 0 & = &
c_{1}\\
\#\{j:v(\sseq(\beta_{j}))=2-2\} & = & \#\{3\} & = & 1 & = & c_{2}.
\end{array}
\]

\end{example}

\begin{proof}[Proof of Lemma \ref{lem:elim_inner_count}]
Recall Lemma \ref{lem:elim_alpha}:%
\begin{equation}
c_{t}=\#\left\{  j: \overline\# \left\{  x:f\left(  x\right)
=0\ \wedge\ x>\beta_{j}\right\}  =m-t\right\}  . \label{eq:elim_alpha}%
\end{equation}
We will eliminate the inner count symbol $\overline\#$ from the right-hand side.
{This follows immediately by applying the Budan-Fourier root counting method \cite{Budan:1807,Fourier:1831} to $f(x)$ (or, equivalently, by applying Descartes' rule of signs to $f(r + \beta_j)$), which is exact because $f$ has only real roots.
Thus, after noting that $f^{(m)} > 0$, we have
\[
  \overline\#\{x\mid f(x)=0\wedge x>\beta_{j}\} = v\left(  \sseq(\beta_{j}) ,+\right).
  \]
  }

\noindent Finally, the claim of the lemma follows immediately by replacing the
expression $\overline\#\left\{  x:f\left(  x\right)  =0\ \wedge\ x>\beta_{j}\right\}  $
with the expression $v\left(  \sseq(\beta_{j}) ,+\right)  $, in the right-hand
side of (\ref{eq:elim_alpha}).
\end{proof}

\bigskip

At this point, we have written each component of $\GEC(F,G)$ as a root counting problem counting eigenvalues of $G$:

\[
c_{t}=\#\left\{  j:v\left(  \sseq(\beta_{j}) ,+\right)
=m-t\right\}. 
\]

\noindent Next, we will rewrite these $m$ root counting problems into matrix form.

\begin{lemma}
  [Rewrite in matrix form]\label{lem:elim_sv}
  {Let $F$ and $G$ be a pair of strongly generic real symmetric matrices. }
  We have
\[
\GEC\left(  F,G\right) = V q
\]
where

\begin{itemize}
\item $V$ is the matrix whose rows are indexed by $t\in\left[  1,\ldots
,m\right]  $ and the columns are indexed lexicographically by $s\in\left\{
-,+\right\}  ^{m}$ with $V_{t,s}=\one_{v(s,+),m-t}$.

\item $q$ is the column vector whose rows are indexed lexicographically by
$s\in\left\{  -,+\right\}  ^{m}$ with
\[
q_{s}=\#\left\{  j:\ \sseq(\beta_{j}) =s\right\}
\]

\end{itemize}
\end{lemma}

\begin{example}
[Running, Section \ref{sec:strong}]\label{ex:prev2} Recall the running example
(Example \ref{ex:strong}), where
$c    = (0,1).$

\noindent From the definition of $V$ in the statement of Lemma \ref{lem:elim_sv}, we have
\[
V=%
\begin{bmatrix}
  1 & 1 & 0 & 0 \\ 0 & 0 & 0 & 1
\end{bmatrix}
.
\]
and $q$ is the vector
\[
q=%
\begin{bmatrix}
  q_{- -} \\
  q_{- +} \\
  q_{+ -} \\
  q_{+ +} \\
\end{bmatrix}
=%
\begin{bmatrix}
\#\{j:\sseq(\beta_{j})=(- -)\}\\
\#\{j:\sseq(\beta_{j})=(- +)\}\\
\#\{j:\sseq(\beta_{j})=(+ -)\}\\
\#\{j:\sseq(\beta_{j})=(+ +)\}\\
\end{bmatrix}
.
\]
In the last example (Example \ref{ex:prev1}),
we calculated that
$\sseq(\beta_{1})    =(+ -)$, 
$\sseq(\beta_{2})    =(+ -)$, and
$\sseq(\beta_{3})    =(++)$.
Hence
\[
q=%
\begin{bmatrix}
\#\{j:\sseq(\beta_{j})=( - -)\}\\
\#\{j:\sseq(\beta_{j})=( - +)\}\\
\#\{j:\sseq(\beta_{j})=( + -)\}\\
\#\{j:\sseq(\beta_{j})=( + +)\}\\
\end{bmatrix}
=%
\begin{bmatrix}
\#\emptyset\\
\#\emptyset\\
\#\{1,2\}\\
\#\{3\}
\end{bmatrix}
=%
\begin{bmatrix}
0\\
0\\
2\\
1
\end{bmatrix}
.
\]
Putting it together, we get
\[
Vq=%
\begin{bmatrix}
  1 & 1 & 0 & 0 \\ 0 & 0 & 0 & 1
\end{bmatrix}%
\begin{bmatrix}
0\\
0\\
2\\
1
\end{bmatrix}
=%
\begin{bmatrix}
0\\
1
\end{bmatrix}
=%
\begin{bmatrix}
c_{1}\\
c_{2}
\end{bmatrix}.
\]

\end{example}

\begin{proof}[Proof of Lemma \ref{lem:elim_sv}]
Recall Lemma \ref{lem:elim_inner_count}:%
\begin{equation}
\label{eq:elim_inner_count}c_{t}=\#\left\{  j:v\left(
\sseq(\beta_{j}) ,+\right)  =m-t\right\}  .
\end{equation}
We will eliminate the sign variation count symbol from $v$ from the right-hand
side. Let $J_t$ denote the set whose entries are counted on the right-hand side of (\ref{eq:elim_inner_count}); that is,
\[
J_{t}=\left\{  j:\ v\left(  \sseq(\beta_{j}) ,+\right)  =m-t\right\}  .
\]
We partition the set $J_{t}\ $according to the sign vector, obtaining
\[
J_{t}=\biguplus\limits_{\substack{s\in\left\{  -,+\right\}  ^{m}\\v\left(
s,+\right)  =m-t}}J_{s}\ \ \ \ \ \text{where }J_{s}=\left\{  j:\ \sseq(\beta
_{j}) =s\right\}  .
\]

\noindent Note that this is indeed a disjoint union, since by the generic
and strongly generic conditions we have~${f^{(k)}(\beta_{j}) \ne0}$ for all $k =
0, \dots, m-1$.
Then we immediately have%
\[
c_{t}=\#J_{t}=\sum_{\substack{s\in\left\{  -,+\right\}  ^{m}\\v\left(
s,+\right)  =m-t}}\#J_{s}%
\]
In matrix form, we can write this as
$c= V q$
where

\begin{itemize}
\item $V$ is the matrix where rows are indexed by $t\in\left[  1,\ldots
,m\right]  $ and the columns are indexed lexicographically by $s\in\left\{
-,+\right\}  ^{m}$ with $V_{t,s}=\one_{v(s+),m-t}$.

\item $q$ is the column vector where rows are indexed lexicographically by
$s\in\left\{  -,+\right\}  ^{m}$ with $q_{s}=\#J_{s}$.
\end{itemize}
Hence we have $\GEC(F,G) = Vq$.
\end{proof}

\begin{lemma}
[Rewrite $q$ in terms of signature]\label{lem:elim_count_N} {Let $F$ and $G$ be a pair of strongly generic real symmetric matrices.} We have
\[
q=H^{-1}\As(F,G)
\]
where $\As(F,G)$ is the vector defined in Definition \ref{def:A}; that is, the vector indexed by $e \in \{0,1\}^m$ where
\[
  (\As(F,G))_e = \sigma(f_e(G))
\]
where
\[
 f_{e}   =\left(  f^{\left(  0\right)  } \right)  ^{e_{0}}\cdots\left(
 f^{\left(  m-1\right)  }\right)  ^{ e_{m-1}}.
\]
\end{lemma}

\begin{example}
[Running, Section \ref{sec:strong}]\label{ex:prev3} In Example \ref{ex:prev2},
we computed
$q = (0,0,2,1)^T$.
From the definition of $H$, we compute
\begin{align*}
H^{-1}  &  =  
\begin{bmatrix}
\frac14 & - \frac14 & -\frac14 & \frac14\\
\noalign{\vspace*{1mm}} \frac14 & \frac14 & -\frac14 & -\frac14\\
\noalign{\vspace*{1mm}} \frac14 & -\frac14 & \frac14 & -\frac14\\
\noalign{\vspace*{1mm}} \frac14 & \frac14 & \frac14 & \frac14
\end{bmatrix}
.
\end{align*}
We then compute
\begin{align*}
f_{00}  &  = \left(  f^{(0)}\right)  ^{0} \left(  f^{(1)}\right)  ^{0} = 1\\
f_{01}  &  = \left(  f^{(0)}\right)  ^{0} \left(  f^{(1)}\right)  ^{1} =
2(x-4)\\
f_{10}  &  = \left(  f^{(0)}\right)  ^{1} \left(  f^{(1)}\right)  ^{0} =
(x-4)^2\\
f_{11}  &  = \left(  f^{(0)}\right)  ^{1} \left(  f^{(1)}\right)  ^{1} =
2(x-4)^3.
\end{align*}
Evaluating those polynomials on the matrix $G$, we compute
\begin{align*}
f_{00}(G)  &  =
             \operatorname*{diag}(1,1,1)
\\
f_{01}(G)  &  =
             \operatorname*{diag}(-4, -4, 8)
\\
f_{10}(G)  &  =
             \operatorname*{diag}(4,4,16)
\\
f_{11}(G)  &  =
             \operatorname*{diag}(-16, -16, 128).
.
\end{align*}
So we have \begin{table}[h]
\centering
\begin{tabular}
[h]{c|c|c|c}%
$e$ & Eigenvalues of $f_{e}(G)$ & \# positive eigenvalues & \# negative
eigenvalues\\\hline
$00$ & $(1,1,1)$ & 3 & 0\\
$01$ & $(-4,-4,8)$ & 1 & 2\\
$10$ & $(4,4,16)$ & 3 & 0\\
$11$ & $(-16,-16,128)$ & 1 & 2
\end{tabular}
\end{table}

\noindent Therefore
\[
\As(F,G)=%
\begin{bmatrix}
\sigma(f_{00}(G))\\
\sigma(f_{01}(G))\\
\sigma(f_{10}(G))\\
\sigma(f_{11}(G))%
\end{bmatrix}
=%
\begin{bmatrix}
3-0\\
1-2\\
3-0\\
1-2
\end{bmatrix}
=%
\begin{bmatrix}
3\\
-1\\
3\\
-1
\end{bmatrix}
.
\]
Finally, we have that
\begin{align*}
H^{-1}\sigma &  =%
\begin{bmatrix}
\frac{1}{4} & -\frac{1}{4} & -\frac{1}{4} & \frac{1}{4}\\
\noalign{\vspace*{1mm}}\frac{1}{4} & \frac{1}{4} & -\frac{1}{4} & -\frac{1}%
{4}\\
\noalign{\vspace*{1mm}}\frac{1}{4} & -\frac{1}{4} & \frac{1}{4} & -\frac{1}%
{4}\\
\noalign{\vspace*{1mm}}\frac{1}{4} & \frac{1}{4} & \frac{1}{4} & \frac{1}{4}%
\end{bmatrix}%
\begin{bmatrix}
3\\
\noalign{\vspace*{1mm}}-1\\
\noalign{\vspace*{1mm}}3\\
\noalign{\vspace*{1mm}}-1
\end{bmatrix}
=%
\begin{bmatrix}
0\\
0\\
2\\
1
\end{bmatrix}
=q.
\end{align*}

\end{example}

\begin{proof}[Proof of Lemma \ref{lem:elim_count_N}]
Recall $q$ from Lemma \ref{lem:elim_sv}: $q$ is the column vector where
rows are indexed lexicographically by $s\in\left\{  -,+\right\}  ^{m}$ with
$q_{s}=\#\left\{  j:\ \sseq(\beta_{j}) =s\right\}.$
We will rewrite {this expression} in terms of the signatures of certain symmetric matrices constructed from $F$~and~$G$.

  From applying the technique proposed by Ben-Or, Kozen, and Reif in~\cite{BKR1986}, we have%
\[
  Hq = \As(F,G)
\]
where

\begin{itemize}
  
\item 
  $H$ is the $2^m \times 2^m$ Hadamard matrix whose rows are indexed by $e \in \{0,1\}^m$ and whose columns are indexed by $s \in \{-, +\}^m$ where
  $$ (H_m)_{e,s} = s_1^{e_1} \cdots s_m ^{e_m}$$

\item $\As$ is defined as in Definition \ref{def:A}; that is, $\As$ is the column vector whose rows are
indexed lexicographically by $e=\left(  e_{0},\ldots,e_{m-1}\right)
\in\left\{  0,1\right\}  ^{m}$ with
\[
(\As)_{e}=\sigma(f_e(G))
\]
\newline where
\[
f_{e}=f^{\left(  0\right)  e_{0}}\cdots f^{\left(  m-1\right)  e_{m-1}},
\]
and again $\sigma(f_e(G))$ is the signature of the matrix $f_e(G)$.
\end{itemize}
{ We note that since $F$ and $G$ are a strongly generic pair, we have that $f_e(\beta_j) \ne 0$ for all $j$; in other words, for all $e$, the characteristic polynomial of $G$ and the polynomial $f_e$ are relatively prime. This is necessary to directly apply the result from Ben-Or, Kozen, and Reif \cite{BKR1986} cited above. In Sections \ref{sec:generic} and \ref{sec:arbitrary}, we show that one does not need this condition.} 

Finally, since $H$ is a Hadamard matrix, it is invertible, so we therefore have
\[
  q = H^{-1} \As(F,G)
\]
and we are done.
\end{proof}


\bigskip\noindent Finally, we have arrived at an expression for $\GEC(F,G)$ which does not contain any references to the eigenvalues $\alpha$ or $\beta$. Hence, we are ready to prove the main result (Theorem \ref{thm:main}) for strongly generic $F$~and~$G$,
namely Lemma~\ref{lem:strong}.

\begin{proof}[Proof of Lemma~\ref{lem:strong}]
Let $F\in\mathbb{R}^{m\times m}$ and $G\in
\mathbb{R}^{n\times n}$ be real symmetric \textbf{strongly generic} matrices. We need to prove the following:
\[
\GEC(F,G)\,\,=\,\,\Cs\,\,\As(F,G).
\]

\noindent For this, we begin by recalling the following lemmas from above.

\begin{enumerate}
\item \label{sg:1} From Lemma \ref{lem:elim_sv}, we have
\[
\GEC\left(  F,G\right) = {V} {q}
\]
where ${V}$ is the matrix where rows are indexed by $t\in\left[
1,\ldots,m\right]  $ and the columns are indexed lexicographically by
$s\in\left\{  -,+\right\}  ^{m}$ with
\[
V_{t,s}=\one_{v(s,+),m-t}.
\]

\item \label{sg:2} From Lemma \ref{lem:elim_count_N}, we have%
\[
q=H^{-1}\As(F,G)
\]
where $\As$ is a column vector whose rows are indexed lexicographically
by $e\in\left\{  0,1\right\}  ^{m}$.
\end{enumerate}

\noindent Together, we finally have
\[
\GEC(F,G)\,\,\underset{\ref{sg:1}}{=}\,\,Vq\,\,\underset{\ref{sg:2}}{=}\,\,V\left( H^{-1} \As(F,G) \right)\,\,\underset{\text{Def. }\ref{def:C}}{=}\,\,\Cs\,\,\As (F,G).
\]
We have proved Lemma~\ref{lem:strong}.
~\end{proof}

\subsection{Proof for generic pair of matrices}

\label{sec:generic}

In this subsection, we will extend the result on strongly generic matrices
(Lemma \ref{lem:strong} from the previous subsection) to a slightly more general result for matrices which are generic but not necessarily strongly generic. 
\begin{lemma}
  \label{lem:generic}
  Let $F$ and $G$ be generic real symmetric matrices (i.e., matrices which do not share eigenvalues).
  Then
\[
  \GEC(F, G) =
 \Cs\,\,\As(F,G).
  \]
\end{lemma}
{ 
  The major difficulty extending Lemma \ref{lem:strong} is that the Ben-Or, Kozen, and Reif root counting algorithm from \cite{BKR1986} that was applied in Lemma \ref{lem:elim_count_N} requires that $f_e(\beta_j)$ being nonzero for all $j$; in other words, that $F$ and $G$ are a strongly generic pair of matrices.
  In this section, we will show that, as a result of the structure of the eigenvalue configuration problem, we can drop this condition and apply the result from Ben-Or, Kozen, and Reif as a black box without modification. 
  }
  To do this, we take the following approach.
  \begin{enumerate}
  \item In Lemma \ref{lem:perturb}, we will show that, given generic $F$~and~$G$, we can ``safely'' perturb $F$ into a new real symmetric matrix $\widehat{F}$, so that $\widehat{F}$ and ${G}$ are strongly generic and the eigenvalue configuration of $F$~and~$G$ equals that of $\widehat{F}$ and ${G}$.
    Since $\widehat{F}$ and ${G}$ are strongly generic, we can apply Lemma \ref{lem:strong}.
    
  \item In Lemma \ref{lem:wrap}, we will show that $\Cs\,\,\As(F,G) = \Cs\,\,\As(\widehat{F}, {G})$. However, before we can do that, we need to prove an intermediate result in Lemma \ref{lem:gamma}, which generalizes Lemma \ref{lem:elim_count_N} from the previous subsection.

    Once this is done, the proof will be complete.
  \end{enumerate}
  Our strategy is illustrated by the following diagram.
\[
\begin{array}{|ccc|}
  \hline
  & &\\
  \GEC(F,G) & \overset{{\color{green}?}}{{\color{green}=\mathrel{\mkern-3mu}=}} & \Cs\,\,\As(F,G) \\[1.5ex]
  \text{Lemma \ref{lem:perturb}}\hspace{1em}\Big\Vert \hspace{6em}        &  & \hspace{6em} \Big\Vert \hspace{1em} \text{Lemma \ref{lem:wrap}} \\[1.5ex]
  \GEC(\widehat{F},{G}) &  =\mathrel{\mkern-3mu}= & \Cs\,\,\As(\widehat{F}, {G})\\[1.5ex]
   & \text{Lemma \ref{lem:strong}} & \\
  &&\\
                                                                            \hline
\end{array}
\]

\bigskip 
First, we begin by showing that, given generic $F$~and~$G$, it is always possible to safely (i.e. without disturbing the eigenvalue configuration) shift the eigenvalues of $F$ slightly to form a new symmetric matrix~$\widehat{F}$ so that $\widehat{F}$ and $G$ are strongly generic.

\begin{lemma}[Safely perturb generic matrices to strongly generic matrices]
\label{lem:perturb} Let $F, G$ be generic matrices. Then
there exists a real symmetric matrix $\widehat{F}$ so that

\begin{itemize}
\item $\widehat{F}$ and ${G}$ are strongly generic and

\item $\GEC(F,G)=\GEC(\widehat{F},{G})$.
\end{itemize}
\end{lemma}

\begin{example}
  [Running, Section \ref{sec:generic}]
  Recall the generic, but not strongly generic, matrices $F$~and~$G$ from Example \ref{ex:generic}:
$F=%
\begin{bmatrix}
1 & 0\\
0 & 5
\end{bmatrix}
\in\mathbb{R}^{2\times2}$ and $G=
\begin{bmatrix}
3 & 0 & 0\\
0 & 3 & 0\\
0 & 0 & 6
\end{bmatrix}
\in\mathbb{R}^{3\times3}$
having eigenvalue configuration
${\GEC(F,G)=(2,1).}$
  Now, let $\varepsilon=0.5$ and set
\[
\widehat{F}=F+\varepsilon I=%
\begin{bmatrix}
1+\epsilon & 0 \\
0 & 5+\varepsilon \\
\end{bmatrix}
=%
\begin{bmatrix}
1.5 & 0 \\
0 & 5.5 \\
\end{bmatrix}
.
\]
The derivatives of the characteristic polynomial $\widehat{f}$ of the new matrix $\widehat{F}$ are
\begin{align*}
  \widehat{f}^{(0)}(x) &= x^2 - 7x + 8.25 \\
  \widehat{f}^{(1)}(x) &= 2x - 7. 
\end{align*}
The arrangement of the eigenvalues of $F$, $G$, and $\widehat{F}$ is shown in
the following diagram.
\[
  \begin{tikzpicture}
    \draw (0, 0) -- (7.5, 0); \fill [red] (1,0) circle(3pt)
    node[label=$\alpha_1$]{}; \fill [red] (5,0) circle(3pt)
    node[label=$\alpha_2$]{}; \fill [green] (1.5,0) circle(3pt)
    node[label=$\widehat \alpha_1$]{}; \fill [green] (5.5,0) circle(3pt)
    node[label=$\widehat \alpha_2$]{}; \fill [blue] (3,0) circle(3pt)
    node[label=$\beta_1$]{}; \fill [blue] (3,0.8) circle(3pt)
    node[label=$\beta_2$]{}; \fill [blue] (6,0) circle(3pt)
    node[label=$\beta_3$]{};
  \end{tikzpicture}
  \]
  
\noindent Note that $\widehat{ F }$ and ${G}$ are strongly generic because
$\widehat{f}^{(0)}$ and $\widehat{f}^{(1)}$ are nonzero at  ${\beta}_{1},{\beta}_{2}$, and ${\beta}_{3}$.
 Further, because $\varepsilon$ was chosen small
enough, no eigenvalue of $\widehat{F}$ \textquotedblleft crossed over\textquotedblright%
\ any $\beta$ and so the eigenvalue configuration of $\widehat{F}$ and ${G}$
remains the same as that of $F$~and~$G$.

\end{example}

\begin{proof}
  [Proof of Lemma \ref{lem:perturb}]
  If $F$ and $G$ are a pair of strongly generic matrices, then clearly $\widehat{F} = F$ satisfies the conclusion of the lemma; thus, for the remainder of this proof, assume without loss of generality that the pair $F$ and $G$ is generic but not strongly generic.

  We will construct $\widehat{F}$ from $F$ by
shifting the eigenvalues of $F$ by a positive real number~$\varepsilon$. Let
$\widehat{F}\ :=\ F+\varepsilon I.$
Obviously $\widehat{F}$ is symmetric and its eigenvalues are simply the
eigenvalues of $F$ plus $\varepsilon$. Now we need to
choose~$\varepsilon$ so that the claims of the lemma hold. Let us choose
$\varepsilon$ as follows. 
\[
\varepsilon\ :=\ \frac{1}{2}\min(A,B,C)
\]
where%
\[%
\begin{array}
[c]{lllll}%
A & = & \min\{|\alpha_{i}-\alpha_{i+1}| & : & i=1,\dots,m-1 \text{ and } \alpha_i \ne \alpha_{i+1}\}\\
B & = & \min\{|\alpha_{i}-\beta_{j}| & : & i=1,\dots,m\text{ and }%
j=1,\dots,n\}\\
C & = & \min\{|\beta_{j}-\gamma| & : & j=1,\dots,n\text{ and }f^{(k)}%
(\gamma)=0\text{ for some $1\le k \le m$ and $\gamma\neq\beta_{j}$}\}
\end{array}
\]
By construction, $A, B, C,$ and therefore $\varepsilon$ are positive real numbers. Let us now check that $\widehat{F}$ and ${G}$ satisfy
the two claims of the lemma, one by one.

\begin{itemize}
\item $\widehat{F}$ and ${G}$ are strongly generic.

\begin{enumerate}
\item $\widehat{F}$ and ${G}$ do not share any eigenvalues, because by
construction each eigenvalue of $\widehat{F}$ is, at minimum, a distance of
$\frac{1}{2}\min(B)>0$ from any eigenvalue of $G$.


\item Each eigenvalue of ${G}$ is, at minimum, a distance of
$\varepsilon>0$ from any root of any derivative of $\widehat{f} = \det(xI_m - \widehat{F})$.
\end{enumerate}

\item $\GEC(F,G)=\GEC(\widehat{F},{G})$.

\begin{enumerate}
\item Since $\varepsilon<\min(A)$ and $\varepsilon<\min(B)$, no eigenvalue of
$\widehat{F}$ crosses over an eigenvalue of $F$ or $G$.

\item Hence, the eigenvalue configuration of $\widehat{F}$ and ${G}$ remains the
same as the eigenvalue configuration of $F$~and~$G$.
\end{enumerate}
\end{itemize}
\end{proof}

Now that we have established that we can perturb generic $F$~and~$G$ to get strongly generic $\widehat{F}$ and ${G}$ without disturbing the eigenvalue configuration, we need to show that
\[
  \Cs\,\,\As(F,G) = \Cs\,\,\As(\widehat{F}, {G}).
  \]
  To do this, we first need to extend some notions from the previous section.
  First, to aid in the computation, we need a key lemma which generalizes Lemma \ref{lem:elim_count_N} from the previous section by finding a simple expression for the quantity $H^{-1}\As(F,G)$.
  We first need to define the notion of the boundary and closure of a sign sequence.

\begin{definition}
\ 

\begin{itemize}
\item The boundary of $s \in\{-, +\}^{m}$ , written as $\partial(s)$, is
defined by
\[
\partial(s) = \{ s^{\prime} \in\{-,0,+\}^{m}: s^{\prime} \text{ and } s \text{
differ by one or more zeros } \}.
\]

\item The closure of $s \in\{-, +\}^{m}$, written as $\cl(s)$, is defined by
\[
\cl(s) = \{s\} \cup\partial s.
\]

\end{itemize}
\end{definition}

\begin{example}
Consider the sign sequence $+ -$. We have \ 

\begin{enumerate}
\item $\partial(+-) = \{+0, 0-, 00 \}$

\item $\cl (+-) = \{+-, +0, 0-, 00 \}$.
\end{enumerate}

\end{example}

\begin{remark}
\ 

\begin{itemize}
\item The notion of the boundary of a sign sequence $s$ can be thought of as
the topological boundary of the set in $\mathbb{R}^{m}$ of points whose
coordinates have the signs given by $s$.

\item Note that the number of elements in $\cl(s)$ where $s \in\{-, +\}^{m}$
is $2^{m}$, because we can choose zero or more entries in $s$ to be zero.
\end{itemize}
\end{remark}

\begin{lemma}
  \label{lem:gamma} Let $F,G$ be an {arbitrary (not necessarily generic)} pair of real symmetric matrices.
Then
\[
H^{-1} \As(F, G) =\sum_{j=1}^{n}\gamma_{j}%
\]
where
\begin{align*}
\gamma_{j}  &  =\frac{1}{\#S_{j}}\sum_{s^{\prime}\in S_{j}}e_{s^{\prime}}\\
S_{j}  &  =\left\{  s\in\{-,+\}^{m}:\sseq(\beta_{j})\in\cl(s)\right\} \\
e_{s^{\prime}}  &  =\text{elementary (standard) unit vector with a 1 in the position of
}s^{\prime}\text{ under lex order with } - < +.
\end{align*}
\end{lemma}

\begin{remark}
  Note that the hypothesis of Lemma \ref{lem:gamma} does not require that $F$ and $G$ be generic.
  In this section, we assume that they are generic, but we will revisit this lemma later in the paper when fully generalizing our result to arbitrary matrices.
\end{remark}

\begin{example}
[Running, Section \ref{sec:generic}]Recall that in the running example we
have
\[
\alpha= (1,5), \qquad\qquad\beta= (3,3,6).
\]
Together with the characteristic polynomial of $f$, the eigenvalues of $F$ and $G$ are shown in the following diagram.
$$
  \begin{tikzpicture}[scale=0.75,domain=0:6]
    \draw[<->,color=red,thick] plot[id=exf] function{0.25*(x-1)*(x-5)}
    node[right] {$f(x)$}; \draw[thick] (-0.5,0) -- (7,0) node[right]
    {$x$};
    \node[color=black,align=left] at (-2.3,-1.25)
    {sign$\left(f^{(0)}, f^{(1)}\right)$}; \fill [red] (1,0)
    circle(2pt) node[] {}; \fill [red] (5,0) circle(2pt) node[] {};
    \fill [blue] (3,0) circle(1.5pt) node[below right] {$\beta_2$}; \fill
    [blue] (3,0.2) circle(1.5pt) node[above right] {$\beta_1$}; \fill
    [blue] (6,0) circle(1.5pt) node[below] {$\beta_3$};

    \draw[dashed] (1,-1) -- (1,1.5); \draw[dashed] (3,-1) -- (3,1.5);
    \draw[dashed] (5,-1) -- (5,1.5);

    \node at (0, -1.2) {$+-$}; \node at (1, -1.2) {$0-$}; \node at (2,
    -1.2) {$--$}; \node at (3, -1.2) {$-0$}; \node at (4, -1.2)
    {$-+$}; \node at (5, -1.2) {$0+$}; \node at (6, -1.2) {$++$};
  \end{tikzpicture}
  $$

Lemma \ref{lem:gamma} should be seen as quantifying the contribution $\gamma_j$ of each
eigenvalue $\beta_{j}$ to the eigenvalue configuration vector. To see this, we
will now compute $\gamma_{j}$ for each $j$ in this example.

\begin{table}[H]
\centering
\begin{tabular}
[h]{c|c|c|c}%
$j$ & $\sseq(\beta_{j})$ & $S_{j}$ & $\gamma_{j}$\\[0.7ex]\hline
1 & $- 0$ & $\{- -, - +\}$ & $%
\frac12 \left(
\begin{bmatrix}
1\\
0\\
0\\
0
\end{bmatrix}
+
\begin{bmatrix}
0\\
1\\
0\\
0
\end{bmatrix}
\right) $\\[6ex]
2 & $-0$ & $\{- -, -+\}$ & $\frac12 \left(
\begin{bmatrix}
1\\
0\\
0\\
0
\end{bmatrix}
+
\begin{bmatrix}
0\\
1\\
0\\
0
\end{bmatrix}
\right)  $\\[6ex]
3 & $+ +$ & $\{+ +\}$ & $%
\begin{bmatrix}
0\\
0\\
0\\
1
\end{bmatrix}
$%
\end{tabular}
\end{table}
\noindent Observe that $S_{3}$ is a singleton set. This
happens because $\sseq(\beta_{3})$ does not contain
zeros, because $f^{(k)}(\beta_{3}) \ne0$ for all $k \in\{1, \dots, m-1\}$. Hence, $\gamma_3$ is simply a unit vector with a 1 in the slot corresponding to the sign sequence of $\beta_3$.

However, $S_{1}$ and $S_2$ are sets with two elements each. This is because $\sseq(\beta
_{1})$ and $\sseq(\beta_2)$ each contain a zero, and so lie on the boundary between the region
$\{x\in\mathbb{R}:\sseq(x)=--\}$ and the region~${\{x\in\mathbb{R}:\sseq(x)=-+\}}$. As a result, the contributions of $\beta_{1}$ and $\beta_2$ to the eigenvalue
configuration vector are split equally between those two regions. 
\end{example}

\begin{remark}
In the following proof, we use bars over certain symbols (e.g. $\overline{H}$)
to denote ``augmented'' versions of those objects used previously in the
paper. These are all defined analogously to their corresponding versions used earlier.
\end{remark}

\begin{proof}
[Proof of Lemma~\ref{lem:gamma}] The proof is long, so we will divide the
proof into a few stages.

\begin{enumerate}
\item 
  First, we apply \cite{BKR1986} as in Lemma \ref{lem:elim_count_N}, with the key difference that we now need to use a larger $3^m \times 3^m$ matrix $\overline{H}$ since the sign sequence of $\beta_j$'s may now contain zeros.
  We will then split up the matrix $\overline{H}$ and find an expression for $\As$.

\item Next, we will find an explicit form for the term $B = H^T \widehat{H}$ appearing in the above
expression.

\item Finally, we rewrite $H^{-1}\As$ in a form which allows us to
quantify the contribution of each eigenvalue $\beta_{j}$ to the overall
eigenvalue configuration vector.
\end{enumerate}

\noindent Now, we elaborate the details for each step.

\begin{enumerate}
  
\item First, we apply \cite{BKR1986} as in Lemma \ref{lem:elim_count_N}, with the key difference that we now need to use a larger $3^m \times 3^m$ matrix $\overline{H}$ since the sign sequence of some $\beta_j$'s may now contain zeros.
  We will then split up the matrix~$\overline{H}$ and find an expression for $\As$ in terms of these pieces of $\overline{H}$.

Let
$$\begin{array}{rll}
(\overline H)_{\overline{e}, \overline{s}}  &  =
\overline{s}_1^{\overline{e}_1} \cdots \overline{s}_m^{\overline{e}_m} \qquad &\text{ where } \overline{e} \in \{0,1,2\}^m \text{ and } \overline{s} \in \{-, 0, +\}^m\\[1.5ex]
\overline{\As}_{\overline e}  &  =
{\sigma ( f_{\overline{e}} (G))} \qquad &\text{ where } \overline{e}    \in\{0,1,2\}^{m}.
\end{array}$$
Note that $H$ (as defined earlier) is a submatrix of $\overline{H}$.

Now, if we apply very similar reasoning as used in the proof of Lemma
\ref{lem:elim_count_N}, we have%

\begin{equation}
\overline{\As} = \overline{H} \overline q. \label{eq:ugly}%
\end{equation}

where
\[%
\begin{array}
[c]{lllll}%
\sigma_{e} & = & \sigma ( f_{e} (G)) & \text{ for } e
\in\{0,1\}^{m} & \\
\widehat{\sigma} & = & \overline\sigma\setminus\sigma &  & \\
\overline q_{s} & = & \# \left\{  j : \sseq(\beta_{j}) = s \right\}  & \text{
for } s \in\{-, 0, +\}^{m} & \\
q_{s} & = & \# \left\{  j : \sseq(\beta_{j}) = s \right\}  & \text{ for } s
\in\{-, +\}^{m} & \\
\widehat{q} & = & \overline q \setminus q. &  &
\end{array}
\]

{The  notation  $\overline \sigma \setminus \sigma $ (similarly for $ \overline q \setminus q$) means that from the column vector $\overline \sigma$, we take away the elements from $\sigma$.}
Note that the key difference between (\ref{eq:ugly}) and the reasoning in
Lemma \ref{lem:elim_count_N} is that now we have~$f_{\overline e}(\beta_{j}) =
0$ for some $j$ and some $\overline{e}$, since $F$~and~$G$ are no longer strongly generic. Hence, we must use the ``augmented'' matrix $\overline{H}
\in\mathbb{N}^{3^{m} \times3^{m}}$, rather than the previously used $H
\in\mathbb{N}^{2^{m} \times2^{m}}$. This hurts the space complexity (i.e., the
number of matrices whose signatures need to be considered) by
increasing the factor from $2^{m}$ to $3^{m}$. We will now address this by
rewriting (\ref{eq:ugly}) in such a way as to reduce the sizes of the matrices involved.

In block matrix form, after appropriately rearranging the rows and columns we
have
\[
\underbrace{\mleft[
\begin{array}{c}
\sigma \\
\hline
\widehat \sigma
\end{array}
\mright]}_{\overline\sigma} = \underbrace{\mleft[
\begin{array}{c | c}
H & \widehat H \\
\hline
\ast    & \ast
\end{array} \mright]}_{\overline H} \underbrace{\mleft[
\begin{array}{c}
q \\
\hline
\widehat q
\end{array}
\mright]}_{\overline q}.
\]

Thus
\begin{align*}
\sigma &  = Hq + \widehat{H} \widehat{q}\\
&  = Hq + HH^{-1} \widehat{H} \widehat{q}\\
&  = H (q + H^{-1} \widehat{H} \widehat{q})\\
&  = H \left(  q + \frac{1}{2^{m}} \underbrace{H^{T} \widehat{H}}_{B}
\widehat{q}\right)  .
\end{align*}
Note that the last equality follows because $H$, being a Hadamard matrix, has inverse $\frac{1}{2^m}H^T$.

Hence, we have found an expression for $\sigma$ in terms of $H, q,
\widehat{q}, $ and $\widehat{H}$.

\item 
Next, we will find an explicit form for the term $B = H^T \widehat{H}$ appearing in the above
expression.
Note that in the following, the term $\widehat{s}$ denotes elements of $\{-,0,+\}^m \setminus \{-, +\}^m$, i.e., sign sequences of length~$m$ with at least one zero.   

Note that by construction, the matrix $B$ can be indexed in its rows by $s \in \{-, +\}^m$ and in its columns by $\widehat{s} \in \{-,0,+\}^m \setminus \{-, +\}^m$.
Thus we have
\begin{alignat*}{2}
  B_{s, \widehat{s}} &= \sum_{e \in\{0,1\}^{m}}^{} (H^{T})_{s, e} \widehat{H}_{e, \widehat{s}} && \text{by expanding the matrix multiplication} \\
                     &= \sum_{e \in\{0,1\}^{m}} \left(  s_{1}^{e_{1}} \cdots s_{m}^{e_{m}}\right) \left(  \widehat{s}_{1}^{e_{1}} \cdots\widehat{s}_{m}^{e_{m}} \right) \hspace{2em} && \text{by the definition of $H$ } \\
                     &=\sum_{e \in\{0,1\}^{m}} \prod_{i=1}^{m} (s_{i} \widehat{s}_{i}) ^{e_{i}} && \text{combining terms by exponents } e_i \\
                     &=  \sum_{e \in\{0,1\}^{m}} \prod_{i=1}^{m} s_{i} \widehat{s}_{i} &&  \text{since } e_i \text{ only has values of 0 or 1} \\
                     &=  \sum_{e \in\{0,1\}^{m}} \prod_{\substack{i=1 \\e_{i} = 1}}^{m} s_{i}\widehat{s}_{i} &&  \text{since the terms in the product are 1 if } e_i = 0\\
                     &=  \sum_{T \subseteq[m]}^{} \prod_{t \in T}^{} s_{t} \widehat{s}_{t} &&  \text{since } \#\{0,1\}^m \text{ equals the number of subsets } T \subseteq [m]\\
                     &=  \sum_{T \subseteq Z_{\widehat{s}}} \prod_{t \in T} \begin{cases}
                       -1 & \text{if } s_{t} \ne \widehat{s}_{t}\\
                       1 & \text{else}%
                     \end{cases}
                       && \text{where } Z_{\widehat{s}} := \{ t \in[m] : \widehat{s}_{t} \ne0 \}\\
&=  \sum_{T \subseteq Z_{\widehat{s}}} (-1)^{\# {\{t \in T : s_{t}  \ne\widehat{s}_{t}\}} } && 
 \text{since only multiplications by } -1 \text{ affect the product.}  \\
\end{alignat*}
Let $U_{T} = \{t \in T : s_{t} \ne\widehat{s}_{t} \}$, and let $U =
U_{Z_{\widehat{s}}}$.

Next, note that for each positive integer $k$, the number of subsets $T
\subseteq Z_{\widehat{s}}$ such that $\#(T \cap U) = k$ is
$\binom{\# U}{k} 2^{\# (Z_{\widehat{s}} \setminus U)}.$
Further, note that the summand terms $(-1)^{\#U_T}$ depend only on the size of the sets~${T \cap U}$.
Hence, by partitioning the sum
$\sum_{T \subseteq Z_{\widehat{s}}} (-1)^{\# U_T}$
by the size of $\#(T \cap U)$, we get
\begin{align*}
B_{s, \widehat{s}}  &  = \sum_{T \subseteq Z_{\widehat{s}}} (-1)^{\#U_{T}}
  = \sum_{k=0}^{\#U} (-1)^{k} \binom{\# U}{k} 2^{\# (Z_{\widehat{s}}\setminus U)}.
\end{align*}
Rearranging, we get
\begin{align*}
B_{s, \widehat{s}} &= 2^{\# (Z_{\widehat{s}} \setminus U)} \sum_{k=0}^{\#U} (-1)^{k}\binom{\#U }{k}\\
&  =
\begin{cases}
2^{\# Z_{\widehat{s}}} & \text{if } U = \emptyset\\
0 & \text{if } U \ne\emptyset.
\end{cases}
\end{align*}
Note that by definition of $\partial s$, we have that $U = \emptyset$ if and
only if $\widehat{s} \in\partial s$. Hence
\[
B_{s, \widehat{s}} =
\begin{cases}
2^{\# Z_{\widehat{s}}} & \text{if } \widehat{s} \in\partial s\\
0 & \text{else } .
\end{cases}
\]
With that, we have found an explicit form for the $B$ matrix.

\item Finally, we rewrite $H^{-1}\As$ in a form which allows us to
quantify the ``contribution'' of each eigenvalue~$\beta_{j}$ to the overall
eigenvalue configuration vector. From the above, for
all $s \in\{-,+\}^{m}$ we have that the~$s$-th component of $H^{-1}\As$ is
\begin{align*}
(H^{-1} \As)_s &  = \left(q + \frac{1}{2^{m}} B \widehat{q}\right)_s 
               = q_s + \frac{1}{2^{m}}\sum_{\widehat{s} \in\partial s} 2^{\# Z_{\widehat{s}}} \widehat{q}_{\widehat{s}}\\
             &  = q_{s} + \sum_{\widehat{s} \in\partial s} \left(  \frac{1}{2}\right)  ^{\# \{0 \in\widehat{s} \}} \widehat{q}_{\widehat{s}} \qquad \text{ since } \# Z_{\widehat{s}}\  -\  m \ = \ \#\{0 \in \widehat{s}\}.
\end{align*}
Noting that $\cl(s) = \partial s \cup s$, we can push the $q_s$ term into the summation and get
\begin{align*}
(H^{-1}\As)_{s} &= \sum_{s^{\prime}\in\cl(s)} \left(  \frac12 \right)  ^{\# \{ 0 \in
s^{\prime}\}} \overline{q}_{s^{\prime}}.
\end{align*}
Note that we switch to the symbol $\overline{q}_{s'}$ to denote the fact that $s'$ could come from the set $\{ -,0,+\}^m$ or $\{-,+\}^m$.
Next, we rewrite $\overline{q}_{s'}$ as a summation of indicators.
\begin{align*}
(H^{-1}\As)_{s}  &= \sum_{s^{\prime}\in\cl(s)} \left(  \frac12 \right)  ^{\# \{ 0 \in
s^{\prime}\}} \sum_{j=1}^{n} \one_{ \sseq(\beta_{j}) , s^{\prime}},
\end{align*}
where here $\one$ denotes the Kronecker delta function. 
Continuing, by rearranging the order of summation we have
\begin{align*}
 (H^{-1}\As)_{s} &  = \sum_{j=1}^{n} \sum_{s^{\prime}\in\cl(s)} \left(  \frac12 \right)
^{\# \{ 0 \in s^{\prime}\}}\one_{ \sseq(\beta_{j}) , s^{\prime}}.
\end{align*}
Now, note that the summand is nonzero if and only if $\sseq(\beta_j)=s'$.
If we combine this with the innermost summation over elements of $\cl(s)$, we can rewrite as
\begin{align*}
  (H^{-1}\As)_{s}&  = \sum_{j=1}^{n}\left(  \frac12 \right)  ^{\# \{ 0 \in\sseq(\beta_{j})
\}}\one_{ \sseq(\beta_{j}) \in \cl(s)}.
\end{align*}
Now, fix $j \in[n]$ and recall that $S_{j} = \{ s \in\{-, +\}^{m} :
\sseq(\beta_{j}) \in\cl(s) \}.$ If $\sseq(\beta_{j})$ has no zeroes, then the
set $S_{j}$ consists of exactly one element; namely, $\sseq(\beta_{j})$
itself. Now suppose that $\sseq(\beta_{j}) = (\dots, 0, \dots)$; i.e., there
is at least one zero. Then for each zero, there are two corresponding elements in $S_{j}$,
which are $(\dots, +, \dots)$ and $(\dots, -, \dots)$. Hence
\[
\# S_{j} = 2^{\# \{0 \in\sseq(\beta_{j}) \}}.
\]
Thus
\begin{align*}
 (H^{-1}\As)_{s}  &= \sum_{j=1}^{n}\left(  \frac12 \right)  ^{\# \{ 0 \in
\sseq(\beta_{j}) \}}\one_{ \sseq(\beta_{j}) \in\cl(s)}\\
&  = \sum_{j=1}^{n} \frac{1}{\# S_{j}} \one_{ \sseq(\beta_{j}) \in
\cl(s)}\\
&  = \sum_{j=1}^{n} \frac{1}{\# S_{j}} \sum_{s^{\prime}\in S_{j}}
\one_{s^{\prime},\sseq(\beta_{j})}.
\end{align*}
Therefore, in full vector form, we have
\[
H^{-1} \As= \sum_{j=1}^{n} \frac{1}{\# S_{j}} \sum_{s^{\prime}\in
S_{j}} e_{s^{\prime}} = \sum_{j=1}^{n} \gamma_{j}.
\]
\end{enumerate}
\end{proof}

{ 
\noindent With that, we are ready to complete the final piece of the proof of Lemma \ref{lem:generic}.
\begin{lemma}
\label{lem:wrap} Let $F, G$ be  generic but not strongly generic, and let
$\widehat{F}$ be symmetric. If we have

\begin{enumerate}
\item $\GEC(F, G) = \GEC(\widehat{F}, {G})$ and

\item $\widehat{F}$ and ${G}$ are strongly generic,
\end{enumerate}

\noindent then
\[
\Cs\,\,\As(F, G) = \Cs\,\,\As(\widehat{F}, {G}).
\]

\end{lemma}

\begin{proof}
  Assume that $F, G$ are generic but not strongly generic, and let $\widehat{F}$ be symmetric.
  Suppose that
\begin{enumerate}
\item $\GEC(F, G) = \GEC(\widehat{F}, {G})$ and

\item $\widehat{F}$ and ${G}$ are strongly generic.
\end{enumerate}
From Lemma~\ref{lem:gamma} , we have 
\[H^{-1}\As(F,G) = \sum_{j=1}^{n}\gamma_{j}\]
and
\[H^{-1}%
\As(\widehat{F},{G})=\displaystyle\sum
_{j=1}^{n}\widehat{\gamma_j}\]
Thus, it suffices to show
\[ V\displaystyle\sum_{j=1}^{n}\gamma_{j}  =  V\displaystyle\sum
_{j=1}^{n}\widehat{\gamma}_{j},
\]
equivalently, to show 
\[
\sum_{j=1}^{n}V\gamma_{j}  =  \displaystyle\sum
_{j=1}^{n}V\widehat{\gamma}_{j}.
\]
Below, we will show a stronger result: for all $j=1,\dots,n$, we have
   
\[V\gamma_{j}=V\widehat{\gamma}_{j}.\] 
Let $j\in\{1,\dots,n\}$ be arbitrary but fixed. Recall that
\[
\gamma_{j}=\frac{1}{\#S_{j}}\sum_{s^{\prime}\in S_{j}}e_{s^{\prime}}%
\]
where
\[
S_{j}=\{s\in\{-,+\}^{m}:\sseq(\beta_{j})\in\cl(s)\},
\]
with $\widehat{\gamma}_{j}$ defined analogously with $\widehat{\beta}_{j}$.
Using the definition of $V$ from the proof of Lemma
\ref{lem:elim_sv}, we then have the following (note that $\gamma_j$ can be substituted with $\widehat{\gamma}_j$):
\begin{alignat*}{2}
(V\gamma_{j})_{t}    &=\sum_{\substack{s\in\left\{  -,+\right\}^{m}\\v\left(  s,+\right)  =m-t}}(\gamma_{j})_{s}\\
                &= \sum_{\substack{s\in\left\{  -,+\right\}  ^{m}\\v\left(  s,+\right)=m-t}}\left(  \frac{1}{\#S_{j}}\sum_{s^{\prime}\in S_{j}}e_{s^{\prime}}\right)  _{s} && \text{by definition of } \gamma_j \\
                &  \sum_{\substack{s\in\left\{  -,+\right\}  ^{m}\\v\left(  s,+\right) =m-t}}\frac{1}{\#S_{j}}\one_{s\in S_{j}} && \text{by evaluating the innermost summation} \\
                &=  \frac{1}{\#S_{j}}\sum_{\substack{s\in\left\{  -,+\right\}  ^{m}\\v\left(s,+\right)  =m-t}}\one_{s\in S_{j}}&& \text{by rearranging} \\
&=  \frac{1}{\#S_{j}}\#\{s\in S_{j}:v(s,+)=m-t\}.
\end{alignat*}
We now make the following claim.

\medskip
\noindent \textbf{Claim: } We have
\begin{equation}
  \label{eq:claim}
 \hspace{5em} V\gamma_j= e_{m - v(\sseq_f(\beta_{j}), +)}   \qquad \text{ and } \qquad V\widehat{\gamma
}_j=e_{m-v(\sseq_{\widehat{f}}({\beta}_{j}),+)}.
\end{equation}

\noindent We will show the proof for the former equality, because the proof for the latter is essentially identical. To prove the claim, it suffices to show that for all $s\in S_{j}$, we have
$v(s,+)=v(\sseq(\beta_{j}),+)$.

\medskip
\noindent Let $s \in S_{j}$. Then we have two cases: either $\sseq(\beta_j) = s$ or not.


\begin{description}
\item[\it Case 1: $\sseq(\beta_{j}) = s$.] Then $S_{j} = \{s\}$, and so
  \begin{align*}
    (V\gamma_j)_t &= \frac{1}{\#S_{j}}\#\{s\in S_{j}:v(s,+)=m-t\} \\
             &= \frac{1}{1}\#\{s\in \{s\} : v(s,+)=m-t\} \\
             &= \one_{ v(\sseq(\beta_j), +) , m-t}.
  \end{align*}
Thus, in matrix form, we have $V\gamma_j= e_{m - v(\sseq(\beta_{j}), +)}$ and the claim follows.
\item[\it Case 2: $\sseq(\beta_{j}) \ne s.$] Then $\sseq(\beta_{j})$ and $s$
differ by one or more zeros. At first glance, it is tempting to use the fact
that for all $y \in\mathbb{R}$, we have that
\[
\# \{x : f(x) = 0 \wedge x > y\} = v(\sseq(y), +),
\]
which was proved as part of Lemma \ref{lem:elim_inner_count} using Descartes'
rule of signs. However, it is not always true that there exists some $y
\in\mathbb{R}$ so that $s = \sseq(y)$. In fact, since $S_{j}$ comprises the
sign sequences of open sets in $\mathbb{R}^{m}$ which share a boundary with
the (not open) set $\{x : \sseq(x) = \sseq(\beta_{j})\}$, there are at most
two sign sequences in $S_{j}$ that are actually realizable. As a result, we
need a different approach. We will instead study the possibilities for
$(\sseq(\beta_{j}), +)$ under our genericity assumptions for $F$.

First, recall that $F$ has exactly $m$ (counting with multiplicity) real eigenvalues.
We then have the following.

{ 
  Suppose $h$ is a polynomial with real roots $(\alpha_i)$ with $h(\beta) \ne 0$ and $h'(\beta) = 0$.
  Since $(h' / h)' = -\sum_{i} (x - \alpha_i)^{-2} < 0$, if we evaluate at $x = \beta$ and use the fact that $h'(\beta) = 0$ we get
  \[
    \frac{h^{(2)}(\beta)}{h(\beta)} = - \sum_{i} \frac{1}{(\beta - \alpha_i)^2} < 0.
  \]

  Suppose now that $f^{(k)}(\beta_j) = 0$ for some $k$.
  We show by induction on $k$ that $f^{(k-1)}(\beta_j)$ and $f^{(k+1)}(\beta_j)$ are nonzero and opposite sign.
  Since $(F,G)$ is generic, we have $f(\beta_j) \ne 0$, so $k \ge 1$.
  If $k = 1$, then $f^{(k-1)}(\beta_j) = f(\beta_j) \ne 0$.
  If $k \ge 2$, then $f^{(k-1)}(\beta_j) = 0$ contradicts the induction hypothesis.

  Now, applying the above to $h = f^{(k-1)}$ (which has only real roots by Rolle's theorem), we have that $f^{(k-1)}(\beta_j)$ and $f^{(k+1)}(\beta_j)$ have opposite sign.
  }

  From the above, we have that any zero in $\sseq(\beta_j)$ must be surrounded by opposing signs; i.e., any zero must appear in the subsequence $+, 0, -$ or $-, 0, +$.
Note that in both of these cases, changing the zero to either a $+$ or a $-$ does not change the sign variation count. 
Hence,
it follows that for all~
$s \in S_{j}$, we have $v(s, +) = v(\sseq(\beta_{j}), +)$.
Then
  \begin{align*}
    (V\gamma_j)_t &= \frac{1}{\#S_{j}}\#\{s\in S_{j}:v(s,+)=m-t\} \\
             &= \frac{1}{\#S_{j}}
               \begin{cases}
                 \#S_j & \text{if }v(\sseq(\beta_{j}), +) = m - t \\
                 0 & \text{otherwise.}
               \end{cases}
    \\
             &= \one_{ v(\sseq(\beta_j), +) , m-t}.
  \end{align*}
Thus, in matrix form, we have $V\gamma_j= e_{m - v(\sseq(\beta_{j}), +)}$ and the claim follows.
\end{description}

From the two cases above, the claim (\ref{eq:claim}) is now proved. 
Thus we have shown that
\[
V\gamma_j=e_{m-v(\sseq_f(\beta_{j}),+)}\qquad\text{ and }\qquad V\widehat{\gamma
}_j=e_{m-v(\sseq_{\widehat{f}}({\beta}_{j}),+)}.
\]

\bigskip

Finally, note that 
\begin{alignat*}{2}
  v(\sseq_f(\beta_{j}),+) &= \#\{x:f(x)=0\wedge x>\beta_{j}\} \qquad && \text{by Lemma \ref{lem:elim_inner_count}}\\
                      &=\#\{x:\widehat{f}(x)=0\wedge x>{\beta}_{j}\} && \text{since } \GEC(F,G) = \GEC(\widehat{F}, G) \\
                      &=v(\sseq_{\widehat{f}}({\beta}_{j}),+) && \text{by Lemma \ref{lem:elim_inner_count}}.
\end{alignat*}
Thus $V\gamma_{j}=V\widehat{\gamma}_{j}$, and so $\Cs\,\,\As(F,G)=\Cs%
\As(\widehat{F},{G})$.
We have proved the lemma.
\end{proof}

\bigskip

\noindent Finally, we are now ready to prove the main result for this section (Lemma
\ref{lem:generic}).

\begin{proof}
[Proof Lemma \ref{lem:generic}] Let $F$~and~$G$ be
generic real symmetric matrices. We need to prove that%
\[
\GEC\left(  F,G\right)   =  \Cs\,\,\As
(F,G).
\]
We consider two cases.

\begin{description}
\item[Case 1:] $F$~and~$G$ are strongly generic.

From Lemma \ref{lem:strong}, we have
\[
\GEC\left(  F,G\right)   = \Cs\,\,\As(F,G)
\]
and so we are done.

\item[Case 2:] $F$~and~$G$ are not strongly generic.
Recall the following.

\begin{enumerate}
\item \label{step:first} From Lemma \ref{lem:perturb}, there exists some real symmetric
matrix $\widehat{F}$ so that $\widehat{F}$ and ${G}$ are strongly generic and
$\GEC(F,G)=\GEC(\widehat{F},{G})$.

\item \label{step:second} From Lemma \ref{lem:strong}, we have
$\GEC(\widehat{F},{G})=\Cs\,\,\As(\widehat{F},{G})$.

\item \label{step:third} From Lemma \ref{lem:wrap}, we have $\Cs\,\,\As(\widehat{F},{G})=\Cs\,\,\As(F,G)$.
\end{enumerate}

Putting these together, we therefore have
\[
\GEC(F,G)\,\,\underset{\ref{step:first}}{=}\,\,\GEC(\widehat{F},{G}%
)\,\,\underset{\ref{step:second}}{=}\,\,\Cs\,\,\As (\widehat{F},{G})%
\,\,\underset{\ref{step:third}}{=}\,\,\Cs\,\,\As (F,G),
\]

\end{description}
\end{proof}

{ 
\subsection{Proof for arbitrary pair of matrices}
\label{sec:arbitrary}
In this subsection, we will prove the main theorem (Theorem~\ref{thm:main}) by generalizing the result from the previous section. 
For this, let us recall the claim of the main theorem: if $F$ and $G$ are real symmetric matrices, then
  \begin{align*}
    \EC(F,G) = \Cs\,\,\As(F,G).
  \end{align*}

  \noindent To prove our theorem for arbitrary real symmetric $F$~and~$G$, we take the following approach.
  \begin{enumerate}
  \item In Lemma \ref{lem:decomp}, we show that taking the eigenvalue configuration of arbitrary $F$~and~$G$ is the sum of the configurations of $F$ and each eigenvalue $\beta_j$ of $G$.
  \item In Lemma \ref{lem:single}, we show that $\EC(F, [\beta_j]) = \Cs\,\,\As(F, [\beta_j])$ for each eigenvalue $\beta_j$ of $G$.
  \item Finally, in Lemma \ref{lem:combine}, we show that the sum $\sum_{j=1}^{n} \As(F, [\beta_j])$ equals $\As(F, G)$. The proof is then complete.
  \end{enumerate}

  \medskip
  \noindent Our strategy is illustrated by the following diagram.
\[
\begin{array}{|ccc|}
  \hline
  & &\\
  \EC(F,G) & \overset{{\color{green}?}}{{\color{green}=\mathrel{\mkern-3mu}=}} & \Cs\,\,\As(F,G) \\[1.5ex]
  \text{Lemma \ref{lem:decomp}}\hspace{1em}\Big\Vert \hspace{6em}        &  & \hspace{6em} \Big\Vert \hspace{1em} \text{Lemma \ref{lem:combine}} \\[1.5ex]
  \sum_{j=1}^{n} \EC(F, [\beta_j]) &  =\mathrel{\mkern-3mu}= & \sum_{j=1}^{n} \Cs \As(F, [\beta_j])\\[1.5ex]
   & \text{Lemma \ref{lem:single}} & \\
  &&\\
                                                                            \hline
\end{array}
\]
First, let us revisit Definition \ref{def:EC}. For real symmetric matrices $F$ and $G$, we have
\[
  \EC(F,G) = \frac{1}{2^m} \sum_{d \in \{-\varepsilon, \varepsilon\}^m} \GEC(F_d, G)
\]
where
\begin{align*}
  \varepsilon &= \text{ ``small enough'' positive number} \\
  F_d &= \text{matrix with eigenvalues $\alpha_i + d_i$ for $i = 1, \dots, m$.}
\end{align*}
First, we show that Definition \ref{def:EC} can always be satisfied. To do this, we explicitly construct $\varepsilon$ and $F_d$.
\begin{lemma}[Explicit construction of $\varepsilon$ from Definition \ref{def:EC}]
  Let $F$ and $G$ be real symmetric matrices.
  Let
\begin{align*}
  \varepsilon &= \frac{1}{2} \min_{i,j}
      \begin{cases}
        1 & \text{ if } \alpha_i = \beta_j \\
        |\alpha_i - \beta_j| & \text{ if } \alpha_i \ne \beta_j
      \end{cases} \\
  F_d &= F + A \diag(d) A^T \qquad \qquad \qquad \qquad \text{ for } d \in \{-\varepsilon, \varepsilon\}^m\\
  A &= \text{ orthonormal matrix such that } F = A \diag(\alpha_1, \dots, \alpha_m) A^T.
\end{align*}
Then Definition \ref{def:EC} is satisfied.
\end{lemma}
\begin{proof}
  First, note that since $F$ is real symmetric, it admits an orthogonal eigendecomposition.
  Hence there exists an orthogonal matrix $A$ so that
  \[
    F = A \diag(\alpha) A^T.
  \]
  Then, for $d \in \{-\varepsilon, \varepsilon\}^m$, we have
  \begin{align*}
    F_d &= F + A \diag(d) A^T \\
        &= A \diag(\alpha) A^T + A \diag(d) A^T \\
        &= A \left( \diag(\alpha) + \diag(d) \right) A^T.
  \end{align*}
  Thus the eigenvalues of $F_d$ are exactly $\alpha_i + d_i$, which satisfies Definition \ref{def:EC}.

  To conclude, let 
  \[
  \varepsilon= \frac{1}{2} \min_{i,j}
      \begin{cases}
        1 & \text{ if } \alpha_i = \beta_j \\
        |\alpha_i - \beta_j| & \text{ if } \alpha_i \ne \beta_j
      \end{cases}.
    \]
    It remains to show that this $\varepsilon$ satisfies Definition \ref{def:EC}.
    By construction, this $\varepsilon$ is strictly less than the minimum distance between $\alpha_i$ and $\beta_j$ over all $i$ and $j$ with $\alpha_i \ne \beta_j$.
    Since the eigenvalues of $F_d$ are exactly $\alpha_i \pm \varepsilon$, this means that no eigenvalue of $F$ travels far enough to cross over an eigenvalue of $G$.
    { Finally, note that if there is no pair $(i,j)$ with $\alpha_i \ne \beta_j$, then any $\varepsilon > 0$ is enough, since our goal is to perturb the eigenvalues of $F$ so that $F_d$ and $G$ share no eigenvalues.
    For simplicity of presentation, we choose $\varepsilon = \frac{1}{2}$, but any positive number suffices.}
    Hence, the condition of Definition \ref{def:EC} is satisfied, and we are done.
\end{proof}
\bigskip
\noindent Having established that Definition \ref{def:EC} can always be satisfied, we can now proceed.
  \begin{lemma}[Decompose $\EC$ into eigenvalues]
    \label{lem:decomp}
    For all real symmetric $F$~and~$G$ we have
    \[
      \EC(F,G) = \sum_{j=1}^{n} \EC(F, [\beta_j]).
    \]
  \end{lemma}

  \begin{example}[Running, Section \ref{sec:arbitrary}] \label{ex:g2}
    Recall from the running example that $F$ has eigenvalues ${(\alpha_1, \alpha_2, \alpha_3) = (0,1,1)}$ and $G$ has eigenvalues $(\beta_1, \beta_2) = (1,2)$.
    We computed in Example \ref{ex:g1} that
    \[
      \EC(F,G) =
      \begin{bmatrix}
        1/4 \\ 1/2 \\ 5/4
      \end{bmatrix}.
     \]
     On the other hand, by applying Definition \ref{def:EC} directly, we have that
     \begin{align*}
       \EC(F, [\beta_1]) &=
                      \begin{bmatrix}
                        1/4 \\ 1/2 \\ 1/4
                      \end{bmatrix} \\
       \EC(F, [\beta_2]) &=
                      \begin{bmatrix}
                        0 \\ 0 \\ 1
                      \end{bmatrix},
     \end{align*}
     and
     \begin{align*}
       \EC(F, [\beta_1]) + \EC(F, [\beta_2]) &= 
                      \begin{bmatrix}
                        1/4 \\ 1/2 \\ 1/4
                      \end{bmatrix} + 
                      \begin{bmatrix}
                        0 \\ 0 \\ 1
                      \end{bmatrix} =
                                     \begin{bmatrix}
                                       1/4 \\ 1/2 \\ 5/4
                                     \end{bmatrix} \\
                                   &= \EC(F,G).
     \end{align*}
     
  \end{example}

  \begin{proof}[Proof of Lemma \ref{lem:decomp}]
    We break the proof into two cases.
    \begin{description}
    \item[Case 1: $F$~and~$G$ are generic.] Recall Definition \ref{def:GEC}:
      \begin{align*}
        \EC(F,G) = (c_1, \dots, c_m) \text{ where  }{c_t = \# \{j : \beta_j \in A_t\}}.
      \end{align*}
      Note that 
      \begin{align*}
        c_t &= \# \{j : \beta_j \in A_t\} \\
            &= \sum_{j=1}^{n}
              \begin{cases}
                1 & \text{ if $\beta_j \in A_t$} \\
                0 & \text{ else},
              \end{cases}
      \end{align*}
      and the summand is exactly the $t$-th component of $\EC(F, [\beta_j])$.

    \item[Case 2: $F$~and~$G$ are not generic.] Recall Definition \ref{def:EC}:
      \begin{align*}
        \EC(F,G) &= \frac{1}{2^m} \sum_{d \in \{-\varepsilon, \varepsilon\}} \GEC(F_d, G).
      \end{align*}
      Since $F_d$ and $G$ are generic, by Case 1 we have 
      \begin{align*}
        \EC(F,G) &= \frac{1}{2^m} \sum_{d \in \{-\varepsilon, \varepsilon\}} \GEC(F_d, G) \\
                &= \frac{1}{2^m} \sum_{d \in \{-\varepsilon, \varepsilon\}} \sum_{j=1}^{n} \GEC(F_d, [\beta_j]) \\
                &=  \sum_{j=1}^{n} \underbrace{\frac{1}{2^m} \sum_{d \in \{-\varepsilon, \varepsilon\}} \GEC(F_d, [\beta_j])}_{\text{$=\EC(F, [\beta_j])$ by Definition \ref{def:EC}}} \\
                &= \sum_{j=1}^{n} \EC(F, [\beta_j]).
      \end{align*}
    \end{description}
  \end{proof}


  \begin{lemma}[$n=1$]
    \label{lem:single}
    Let $F$ be an arbitrary real symmetric matrix and let $\beta \in \mathbb{R}$.
    Then
    \[
      \EC(F,[\beta]) = \Cs\,\,\As(F,[\beta]).
      \]
  \end{lemma}
  \begin{example}[Running, Section \ref{sec:arbitrary}]
    Let $F = \diag(0,1,1)$ and let $\beta = 1$.
    In Example \ref{ex:g2}, we found via Definition \ref{def:EC} that
    \begin{align*}
      \EC(F, [\beta]) &=
                   \begin{bmatrix}
                     1/4 \\ 1/2 \\ 1/4
                   \end{bmatrix}.
    \end{align*}
    On the other hand, using Definitions \ref{def:C} and \ref{def:A} directly, we have
    \begin{align*}
      \Cs \As(F, [\beta]) &= VH^{-1}
                        \begin{bmatrix}
                          \sigma(f_{000}(G)) \\
                          \sigma(f_{001}(G)) \\
                          \sigma(f_{010}(G)) \\
                          \sigma(f_{011}(G)) \\
                          \sigma(f_{100}(G)) \\
                          \sigma(f_{101}(G)) \\
                          \sigma(f_{110}(G)) \\
                          \sigma(f_{111}(G)) \\
                        \end{bmatrix} \\
      &=
                        \left[\begin{array}{cccccccc}
0 & 0 & 0 & 0 & 1 & 1 & 1 & 0 
\\
 1 & 1 & 0 & 1 & 0 & 0 & 0 & 0 
\\
 0 & 0 & 0 & 0 & 0 & 0 & 0 & 1 
\end{array}\right]
\left[\begin{array}{cccccccc}
1 & 1 & 1 & 1 & 1 & 1 & 1 & 1 
\\
 -1 & 1 & -1 & 1 & -1 & 1 & -1 & 1 
\\
 -1 & -1 & 1 & 1 & -1 & -1 & 1 & 1 
\\
 1 & -1 & -1 & 1 & 1 & -1 & -1 & 1 
\\
 -1 & -1 & -1 & -1 & 1 & 1 & 1 & 1 
\\
 1 & -1 & 1 & -1 & -1 & 1 & -1 & 1 
\\
 1 & 1 & -1 & -1 & -1 & -1 & 1 & 1 
\\
 -1 & 1 & 1 & -1 & 1 & -1 & -1 & 1 
\end{array}\right]^{-1}
                          \left[\begin{array}{c}
1 
\\
 1 
\\
 0 
\\
 0 
\\
 0 
\\
 0 
\\
 0 
\\
 0 
                          \end{array}\right] 
                        =
                          \begin{bmatrix}
                            1/4 \\ 1/2 \\ 1/4
                          \end{bmatrix}.
    \end{align*} 
  \end{example}
  \begin{proof}[Proof of Lemma \ref{lem:single}]
    Note that if $F$~and~$G$ are generic (i.e. $\beta$ is not an eigenvalue of $F$), then Lemma \ref{lem:generic} applies and this lemma immediately follows as a special case.

    Thus, for the remainder of this proof, suppose that $F$~and~$G$ are not generic.
    In other words, assume that the eigenvalues of $F$ are
    \[
      \alpha_1 \le  \cdots \le \alpha_l < \alpha_{l+1} = \cdots = \alpha_{l + \mu} < \alpha_{l + \mu + 1} \le \cdots \le \alpha_m,
      \]
      where $\alpha_{l+1} = \cdots = \alpha_{l+\mu} = \beta$.
      That is, suppose that $\beta$ is an eigenvalue of $F$ of multiplicity $\mu$, starting at index $l+1$ in the list of eigenvalues of $F$ (repeated with multiplicity).
      \medskip
      
      \noindent To prove the lemma, we will show that 
      \begin{equation}
        \EC(F, [\beta]) = \Cs\,\,\As(F, [\beta]) = 
        \frac{1}{2^\mu} \sum_{k=0}^{\mu} \binom{\mu}{k} e_{l + \mu - k}. \label{eq:goal}
      \end{equation}
      We will first simplify the left-hand side; i.e., we will simplify
      $  \EC(F, [\beta])$
        by opening up Definition \ref{def:EC}.
        We have
    \begin{alignat*}{2}
      \EC(F,[\beta]) &= \frac{1}{2^m} \sum_{d \in \{-\varepsilon, \varepsilon\}^m} \GEC(F_d, G) && \text{by Definition \ref{def:EC}}\\
              &= \frac{1}{2^m} \sum_{d \in \{-\varepsilon, \varepsilon\}^m} \Cs\,\,\As(F_d, [\beta])  \qquad && \text{by Lemma \ref{lem:generic}, since $F_d$ and $[\beta]$ are generic}\\
              &= \frac{1}{2^m} \sum_{d \in \{-\varepsilon, \varepsilon\}^m} V \sum_{j=1}^{1}\gamma_{d,j}  && \text{by Lemma \ref{lem:gamma}}\\
              &= \frac{1}{2^m} \sum_{d \in \{-\varepsilon, \varepsilon\}^m} e_{m - v(\sseq_{f_d}(\beta), +)}  &&  \text{by Claim \ref{eq:claim} in Lemma \ref{lem:wrap}, with $f_d = \det(xI_m-F_d).$}
    \end{alignat*}

    \medskip

    \noindent Next, we will examine the term $v(\sseq_{f_d}(\beta), +)$. Recall that by construction, the eigenvalues of~$F_d$ are~$\alpha_i + d_i$ for~$i \in \{1, \dots, m\}$.
    Recall from Lemma \ref{lem:elim_inner_count} that
      $v(\sseq_{f_d}(\beta), +) = \# \{ i : \alpha_i + d_i > \beta\}.$
    Note further that~${\alpha_i + d_i > \beta}$ if and only if either \begin{enumerate}
    \item $\alpha_i > \beta$,
    \item 
     or $\alpha_i = \beta$ and $d_i > 0$.
    \end{enumerate}
    This means that
    \begin{alignat*}{2}
      \# \{ i : \alpha_i + d_i > \beta \}
      &= \# \{ i : \alpha_i > \beta \} &&+ \,\, \# \{ i : d_i > 0 \,\,\wedge\,\,  \alpha_i = \beta\}  \nonumber\\
      &= m - (l + \mu) &&+ \,\,\# \{ i : d_i > 0 \,\,\wedge\,\,  \alpha_i = \beta\}  \nonumber\\
      &= m - l - \mu &&+ \,\,\# \{ i : d_i > 0 \,\,\wedge\,\,  i \in \{l +1, \dots, l+\mu\} \}.
    \end{alignat*}
    Thus, continuing from above we have
    \begin{align*}
      \EC(F,[\beta]) &= \frac{1}{2^m}  \sum_{d \in \{ -\varepsilon, \varepsilon\}^m} e_{m - v(\sseq_{f_d}(\beta), +)} \\
      &= \frac{1}{2^m} \sum_{d \in \{ -\varepsilon, \varepsilon\}^m} e_{m - (m - l - \mu + \# \{ i : d_i > 0 \,\,\wedge\,\,  i \in \{l +1, \dots, l+\mu\} \})} \\
      &= \frac{1}{2^m} \sum_{d \in \{ -\varepsilon, \varepsilon\}^m} e_{l + \mu - \# \{ i : d_i > 0 \,\,\wedge\,\,  i \in \{l +1, \dots, l+\mu\} \}}.
    \end{align*}
    Next, we will partition the sum by the values of~$\# \{ i : d_i > 0 \,\,\wedge\,\,  i \in \{l +1, \dots, l+\mu\} \}$.
    The possible values for this term are $0$ through $\mu$ inclusive.
    Further, for fixed $k \in \{0, \dots, \mu\}$, the number of $d \in \{-\varepsilon, \varepsilon\}^m$ so that~$\# \{ i : d_i > 0 \,\,\wedge\,\,  i \in \{l +1, \dots, l+\mu\} \} = k$ is exactly $2^{m - \mu}\binom{\mu}{k}$ (i.e. by choosing the arbitrary $m - \mu$ entries, then independently choosing $k$ of the $\mu$ remaining positive entries).
    Therefore
    \begin{align*}
      \EC(F, [\beta]) 
      &= \frac{1}{2^m} \sum_{d \in \{ -\varepsilon, \varepsilon\}^m} e_{l + \mu - \# \{ i : d_i > 0 \,\,\wedge\,\,  i \in \{l +1, \dots, l+\mu\} \}}  \\
      &= \frac{1}{2^m} \sum_{k=0}^{\mu} 2^{m - \mu} \binom{\mu}{k} e_{l + \mu - k}  \\
      &= \frac{1}{2^\mu} \sum_{k=0}^{\mu} \binom{\mu}{k} e_{l + \mu - k}. 
    \end{align*}
    This is the same as (\ref{eq:goal}), so we are done with the left-hand side.

    \medskip
    \noindent Now, we will simplify the right-hand side by rewriting the expression
     $ \Cs\,\,\As(F, [\beta]).$
    Recall that by Lemma~\ref{lem:gamma} we have
    \begin{align*}
      \Cs\,\,\As(F,[\beta]) &= V \gamma
    \end{align*}
    where
    \[
      \gamma = \frac{1}{\# S} \sum_{s \in S} e_{s} 
    \]
    where $S = \left\{  s\in\{-,+\}^{m}:\sseq(\beta)\in\cl(s)\right\}$. 
    Hence we have
    \begin{align*}
                   \Cs\,\,\As(F,[\beta]) &= V\gamma =  V \left( \frac{1}{\# S} \sum_{s \in S} e_{s} \right) \\
                      &= \frac{1}{\#S}
                        \begin{bmatrix}
                          \# \{s \in S: v(s, +) = m - 1 \} \\
                          \# \{s \in S: v(s, +) = m - 2 \} \\
                          \vdots \\
                          \# \{s \in S: v(s, +) = m - m \}
                        \end{bmatrix}.
    \end{align*}
    Recall the proof of Lemma \ref{lem:gamma} for more detailed reasoning on the last step.

    \medskip
    \noindent Next, fix $t \in \{1, \dots, m\}$. We will consider the $t$-th element of the above vector; that is,
    \begin{align}
      \frac{1}{\#S} \{s \in S: v(s, +) = m - t\}. \label{eq:rhs}
    \end{align}
    Recall the definition of $S$ from Lemma \ref{lem:gamma}:
    \[
      S = \{s \in \{-,+\}^m : \sseq_f(\beta) \in \cl(s)\}.
      \]
      Note that since $\beta$ is a root of $f$ with multiplicity $\mu$, we have that
      \[
        \sseq_f(\beta) = (\underbrace{0, \dots, 0}_{\text{$\mu$ zeros}}, \omega),
        \]
        where $\omega \in \{-, 0, +\}^{m - \mu}$.
        Note that $\omega_1 \ne 0$, because by assumption $\beta$ has multiplicity exactly $\mu$.
        Let $z$ denote the number of zeros in $\omega$.
        With this, we have that
        \begin{align*}
          S &= \{ s \in \{-, +\}^{m} : \text{$s$ and $\sseq_f(\beta)$ differ by one or more zeros} \}.
        \end{align*}
        In other words, the set $S$ is the set of sign sequences of length $m$ which can be obtained from $\sseq_f(\beta)$ by replacing each $0$ with either $+$ or $-$.
        Since there are $\mu + z$ zeros in $\sseq_f(\beta)$, this means that $S$ has $2^{\mu + z}$ elements.

        Further, by the same reasoning used in the proof of Lemma \ref{lem:wrap}, {we have that any zeros in $\omega$ appear as $+, 0, -$ or $-, 0, +$; in particular, this means $\omega$ does not have consecutive zeros.}
        Now, let $s \in S$. We will now rewrite $v(s, +) = v(s_1, \dots, s_m, +)$. Note
\begin{alignat*}{2}
  v(s_1, \dots, s_m, +) &= v(s_1, \dots, s_\mu, s_{\mu +1}) + v(s_{\mu + 1}, \dots, s_m, +) \qquad \\
                        &= v(s_1, \dots, s_{\mu}, \omega_1) + v(\omega_1, \dots, s_m, +) && \text{since $\omega_1 \ne 0$, so $s_{\mu+1} = \omega_1$}\\
                        &= v(s_1, \dots, s_\mu, \omega_1) + v(\omega, +) &&\text{by the above property of $\omega$} \\
  &= v(s_1, \dots, s_\mu, \omega_1) + v(\sseq(\beta), +)  && \text{since $v(\omega, +) = v(\sseq(\beta), +)$}\\
  &= v(s_1, \dots, s_\mu, \omega_1) + \# \{i : \alpha_i > \beta \} && \text{since $v(\sseq(\beta), +) = \# \{i : \alpha_i > \beta \}$ by Lemma \ref{lem:elim_inner_count}} \\
  &= v(s_1, \dots, s_\mu, \omega_1) + m - l - \mu &&  \text{by the hypotheses of the current lemma.}
\end{alignat*}
        Now, substituting the above into (\ref{eq:rhs}) gives
        \begin{align*}
          \frac{1}{\#S} \#\{s \in S: v(s, +) = m - t\}
          &= \frac{1}{2^{\mu + z}} \# \{ s \in S: v(s, +) = m - t \} \\
          &= \frac{1}{2^{\mu + z}} \# \{ s \in S:v(s_1, \dots, s_\mu, \omega_1) + m - l - \mu = m - t \} \\
          &= \frac{1}{2^{\mu + z}} \# \{ s \in S:v(s_1, \dots, s_\mu, \omega_1) = l + \mu - t \} 
        \end{align*}
        Note that in counting $\# \{ s \in S : v(s_1, \dots, s_\mu, \omega_1) = l + \mu - t\}$, we choose the first~$\mu$ entries so that there are exactly~$l + \mu - t$ sign variations, and we arbitrarily choose $z$ zeros from the $\omega$ part.
        Hence
        \begin{align*}
          \frac{1}{\#S} \#\{s \in S: v(s, +) = m - t\}
          &= \frac{1}{2^{\mu + z}} 2^z \binom{\mu}{l + \mu - t}\\
          &= \frac{1}{2^\mu} \binom{\mu}{l + \mu - t}. 
        \end{align*}
        Now, we return to the full vector form for $V\gamma$. Note that in the below, we use the convention that~{$m$-dimensional} standard unit vectors $e_i$ are zero when $i \not \in \{1, \dots, m\}$ and binomial coefficients are likewise zero when the denominator is out of range.
        With this, we have
        \begin{align*}
          \Cs\,\,\As(F, [\beta]) &= V\gamma \\
                           &= \frac{1}{2^\mu}
                             \sum_{t}^{} \binom{\mu}{l + \mu - t} e_t \\
                           &= \frac{1}{2^\mu}
                             \sum_{k}^{} \binom{\mu}{k} e_{l + \mu - k} \qquad \text{by reindexing with $k=l+\mu-t$}\\
                           &= \frac{1}{2^\mu} \sum_{k=0}^{\mu} \binom{\mu}{k} e_{l + \mu - k} \qquad \text{since the summand is zero when $k \not \in \{1, \dots, m\}$.}
        \end{align*}
  This is exactly the same as (\ref{eq:goal}); as a result, we have shown that
      \[\EC(F,[\beta])= \frac{1}{2^\mu} \sum_{k=0}^{\mu} \binom{\mu}{k} e_{l + \mu - k}= \Cs\,\,\As(F,[\beta])\]
      and we are done.
  \end{proof}
\noindent  In the next lemma, we show that $\sum_{j=1}^{n} \Cs\,\,\As(F, [\beta_j])$ equals $\Cs\,\,\As(F, G)$.
  \begin{lemma}
    \label{lem:combine}
    Let $F$~and~$G$ be arbitrary real symmetric matrices. Let $\beta_1, \dots, \beta_n$ be the eigenvalues of $G$, including multiple eigenvalues.
    Then
    \[
      \sum_{j=1}^{n} \Cs\,\,\As(F, [\beta_j]) = \Cs\,\,\As(F,G).
      \]
  \end{lemma}

  \begin{proof}[Proof of Lemma \ref{lem:combine}]
    Let $F$~and~$G$ be arbitrary and let $\beta_1, \dots, \beta_n$ be the eigenvalues of $G$.
    Note
    \begin{align*}
      \sum_{j=1}^{n} \Cs \As(F, [\beta_j]) &= \Cs \sum_{j=1}^{n} \As(F, [\beta_j]).
    \end{align*}
    Hence it suffices to show that
    \[
      \sum_{j=1}^{n} \As(F, [\beta_j]) = \As(F, G).
      \]
      Recall that for any two matrices $F$~and~$G$, the column vector $\As(F,G)$ is defined as the vector indexed by~$e \in \{0,1\}^m$ where
      \begin{align*}
        (\As(F,G))_e &= \sigma (f_e(G)) \\
                   &= \# \{j : f_e(\beta_j) > 0 \} - \# \{j : f_e(\beta_j) < 0\}.
      \end{align*}
      Hence
      \[
        (\As(F, [\beta_j]))_e = 
                     \begin{cases}
                       1 & \text{ if $f_e(\beta_j) > 0$} \\
                       -1 & \text{ if $f_e(\beta_j) < 0$}.
                     \end{cases} 
      \]
      Note further that
      \begin{align*}
        (\As(F,G))_e &= \# \{j : f_e(\beta_j) > 0 \} - \# \{j : f_e(\beta_j) < 0\} \\
                   &= \sum_{j=1}^{n}
                     \begin{cases}
                       1 & \text{ if $f_e(\beta_j) > 0$} \\
                       0 & \text{ else}
                     \end{cases} - 
 \sum_{j=1}^{n}
                     \begin{cases}
                       1 & \text{ if $f_e(\beta_j) < 0$} \\
                       0 & \text{ else}
                     \end{cases}  \\
        &= \sum_{j=1}^{n}
                     \begin{cases}
                       1 & \text{ if $f_e(\beta_j) > 0$} \\
                       -1 & \text{ if $f_e(\beta_j) < 0$}
                     \end{cases}  \\
                   &= \sum_{j=1}^{n} \sigma \left( f_e\left([\beta_j]\right) \right) 
                   = \sum_{j=1}^{n} (\As(F, [\beta_j]))_e.
      \end{align*}
      Since this holds for all $e \in \{0,1\}^m$, we have
      \[
        \sum_{j=1}^{n} \As(F, [\beta_j]) = \As(F,G),
        \]
        and we are done.
  \end{proof}
  \medskip
  \noindent With that, we are now ready to prove our main theorem (Theorem \ref{thm:main}) for arbitrary real symmetric matrices $F$~and~$G$.
  \begin{proof}[Proof of Main Theorem (Theorem \ref{thm:main})]
    Let $F$~and~$G$ be arbitrary real symmetric matrices.
    Then
\[
\EC(F,G)\,\,\underset{\text{Lem. }\ref{lem:decomp}}{=}\,\,\sum_{j=1}^{n} \EC(F, [\beta_j])\,\,\underset{\text{Lem. }\ref{lem:single}}{=}\,\,\sum_{j=1}^{n} \Cs\,\,\As(F, [\beta_j])
\,\,\underset{\text{Lem. }\ref{lem:combine}}{=}\,\,\Cs\,\,\As (F,G).
\]
    
  \end{proof}
 } 

 { 
 \section{Application: Rank-$r$ updates}
 \label{sec:update}
 In this section, we briefly discuss one potential application for Theorem \ref{thm:main} to the characterization of eigenvalue configurations which can arise via rank updates of real symmetric matrices.

 \begin{definition}
   Let $F, G \in \mathbb{R}^{m \times m}$ be symmetric with respective eigenvalues ${\alpha_1 \le \cdots \le \alpha_m}$ and ${\beta_1 \le \cdots \le \beta_m}$.
   We say $G$ is \textbf{spectrally reachable by a rank-$r$ update} from $F$ if there exists a positive semidefinite~$E \in \mathbb{R}^{m \times m}$ with rank $E \le r$ such that
     \[
       \operatorname*{spec}(G) = \operatorname*{spec}(F+E),
       \]
       where $\operatorname*{spec}$ denotes the multiset of eigenvalues of a matrix.
 \end{definition}
 \noindent This leads us to the following problem:
\begin{problem}
\ 
\begin{enumerate}
\item[In:\ ] $F, G \in \mathbb{R}[p]^{m \times m}$ symmetric matrices
where $p$ is a finite set of parameters.

$r \in \mathbb{Z}$, a positive integer
\item[Out:] Condition on $p$ equivalent to whether $G$ is spectrally reachable by a rank-$r$ update from $F$
\end{enumerate}

\end{problem}
This problem is addressed via Weyl's inequalities and their converse (\cite{buch2006eigenvalues}): that is, $G$ is reachable via rank-$r$ update from $F$ if and only if
\begin{equation}
  \alpha_i \le \beta_i \le \alpha_{i + r} \label{eq:interlacing}
\end{equation}
for all $i \in \{1, \dots, m\}$.
 For parametric $F$ and $G$, Weyl's inequalities give only quantified semialgebraic conditions, which is not suitable for symbolic computation in terms of the parameters.
 However, to our knowledge, there is no such classical result in terms of the \textit{entries} of the matrices (i.e., a quantifier-free condition).
 We can use Theorem~\ref{thm:main} to give such a quantifier-free condition.

 \begin{corollary}[Reachability of a rank-$r$ update] \label{cor:rank}
   Let $F, G \in \mathbb{R}^{m \times m}$ be symmetric with respective characteristic polynomials $f$ and $g$ and let $(F,G)$ be generic (i.e. the resultant of $f$ and $g$ is nonzero). Let $r \ge 1$ be a positive integer.
   Let $w^{(k)} \in \mathbb{Q}^{2^m}$ be the sum of the first $k$ rows of $\Cs$ and let $T_k = \langle w^{(k)}, \As(F, G) \rangle$.
   Then
   $G$ is spectrally reachable by a rank-$r$ update from $F$ if and only if  we have
   \begin{equation}
     \label{eq:bounds}
     T_m = m \,\,\,\, \wedge \,\,\,\, k + 1 - r \le T_k \le k \text{ for all $k \in \{1, \dots, m-1\}$}.
   \end{equation}
 \end{corollary}
 Note that $T_k$ is a rational linear combination of signatures $\sigma(f_e(G))$, which can be expressed solely in terms of the parameters $p$. Therefore the condition on $T_k$ is quantifier-free.

 \begin{proof}
   By Theorem \ref{thm:main} we have that 
   \[
     EC(F,G) =
     \begin{bmatrix}
       c_1 \\ \vdots \\ c_m
     \end{bmatrix} = \Cs \As (F,G).
   \]
   Hence, it follows that
   \[
     T_k = \langle w^{(k)}, \As(F,G) \rangle = c_1 + \dots + c_k.
   \]
   Note that this means that $T_k = \# \{j : \alpha_1 < \beta_j < \alpha_{k+1} \}$, and so when $T_m = m$ we have
   \begin{equation}
     \label{eq:ineq}
     \beta_j < \alpha_{k+1} \,\, \iff \,\, j \le T_k.
   \end{equation}

   \noindent First, suppose that $G$ is spectrally rank-$r$ reachable from $F$. Then (\ref{eq:interlacing}) holds for all $j$, and since $(F,G)$ is a generic pair, all inequalities are strict.
   Since $\alpha_1 < \beta_1$, we have that
   \begin{align*}
     T_m &= c_1 + \dots + c_m   \\
         &= \# \{ j : \alpha_1 < \beta_j < \alpha_2 \} + \cdots + \# \{j : \alpha_{m} < \beta_j < \infty \} \\
         &= \# \{ j : \alpha_1 < \beta_j\} \\
         &= m.
   \end{align*}
   Fix $k \in \{1, \dots, m-1\}$.
   \begin{description}
   \item[Upper bound ($T_k \le k$):] 
   We have
   $   \alpha_{T_k} < \beta_{T_k} < \alpha_{k + 1}$, which implies that $T_k \le k$.
 \item[Lower bound ($k+1 -r \le T_k$):] 
   Note that 
   \[
     \beta_{k + 1 - r} < \alpha_{(k + 1 - r) + r} = \alpha_{k+1},
    \]
    and therefore $k + 1 - r \le T_k$ by (\ref{eq:ineq}).
   \end{description}
    Therefore
    \[
     T_m = m \,\,\,\, \wedge \,\,\,\, k + 1 - r \le T_k \le k.
      \]

      Conversely, suppose that (\ref{eq:bounds}) holds for all $k \in \{1, \dots, m-1\}$. Fix $j \in \{1, \dots, m\}$.
      Since $T_m = m$, we have that $\alpha_1 < \beta_j$. 
      \begin{description}
      \item[Upper bound ($\beta_j < \alpha_{j + r}$):] 
      Suppose that $j+r \le m$. Since $r \ge 1$, we have that
      $     1 \le j + r - 1 \le m - 1.$
      This means that
      $T_k \ge k+1-r = j$, so by (\ref{eq:ineq}) we have $\beta_j < \alpha_{j+r}$.
      If $j + r > m$, then $\alpha_{j+r} = \infty$ and we are done.
    \item[Lower bound ($\alpha_j < \beta_j$):] 
      Note that by (\ref{eq:bounds}) we have $T_{j-1} \le j-1 < j$, and so by (\ref{eq:ineq}) we have $\alpha_j \le \beta_j$. Since $\beta_j \ne \alpha_j$, this inequality is strict.
      \end{description}
      Putting these bounds together, we get that
      \[
        \alpha_j < \beta_j < \alpha_{j+r}
        \]
        and we are done.
 \end{proof}

\begin{example}[Rank-$r$ reachability with one parameter]
Let
\[
F = \begin{bmatrix} 0 & 0 \\ 0 & 2 \end{bmatrix},
\qquad
G = \begin{bmatrix} p & 1 \\ 1 & p \end{bmatrix},
\]
where $p$ is a parameter such that $F$ and $G$ are a generic pair. We will use Corollary \ref{cor:rank} to determine a quantifier-free condition on $p$ so that $G$ is reachable from $F$ via a rank-$r$ update.
We first compute:
\begin{description}
\item[$\As(F,G)$:] We have
  \begin{align*}
    f_{00}(G) &=
                \begin{bmatrix}
                  1 & 0 \\ 0 & 1
                \end{bmatrix} \\
    f_{01}(G) &=
                \begin{bmatrix}
                  2p - 2 & 2 \\ 2 & 2p-2
                \end{bmatrix}\\
    f_{10}(G) &= 
                \begin{bmatrix}
                  (p-1)^2 & 2p-2 \\
                  2p-2 & (p-1)^2
                \end{bmatrix} \\
    f_{11}(G) &=
                \begin{bmatrix}
                  2(p-1)(p^2-2p+3) & 6(p-1)^2 \\ 6(p-1)^2 & 2(p-1)(p^2 - 2p + 3)
                \end{bmatrix}.
  \end{align*}
\item[$\Cs$:] From the $\Cs$ computed in Example \ref{ex:C}, we have that
  \begin{align*}
    w^{(1)} =
    \begin{bmatrix}
      \frac12 & 0 & - \frac12 & 0
    \end{bmatrix} , \qquad \qquad
    w^{(2)} =
    \begin{bmatrix}
      \frac34 & \frac14 & - \frac14 & \frac14 
    \end{bmatrix}.
  \end{align*}
  Let $\sigma_e := \sigma(f_e(G))$.
  We therefore have
  \begin{align*}
    T_1 &= 1 - \frac12 \sigma_{10} \\
    T_2 &= \frac{3}{2} + \frac14 (\sigma_{01} - \sigma_{10} + \sigma_{11}).
  \end{align*}
\end{description}
Applying Corollary \ref{cor:rank}, we have that $G$ is reachable from $F$ via a rank-$r$ update if and only if
\begin{alignat*}{3}
  &T_2 = 2 \qquad &&\scalebox{2}{$\wedge$}&& \qquad 2 - r \le T_1 \le 1\\
  \iff& \frac{3}{2} + \frac14 (\sigma_{01} - \sigma_{10} + \sigma_{11}) = 2 \qquad &&\scalebox{2}{$\wedge$}&& \qquad 2-r \le 1 - \frac12 \sigma_{10} \le 1.
\end{alignat*}

Finally, we now compute the specific values of $p$ for which the above condition is satisfied.
Note that for a $2 \times 2$ matrix $M$, we have that
\begin{align*}
  \sigma(M) &=
         \begin{cases}
           2 \operatorname*{sign} (\operatorname*{tr} M) & \text{if $\det M > 0$} \\
            \operatorname*{sign} (\operatorname*{tr} M) & \text{if $\det M = 0$} \\
           0 & \text{if $\det M < 0$}. 
         \end{cases}
\end{align*}
Note also that
\begin{align*}
  \det f(G) &= (p-1)^2(p+1)(p-3)
\end{align*}
and so $(F,G)$ is generic if and only if $p \ne \pm 1$ and $p \ne 3$.
In the table below, we compute signatures, values for $T_1$ and $T_2$, and reachability for each valid set of values for $p$.
\[
\begin{array}{lcccccc}
  \hline
p & (\sigma_{00}, \sigma_{01}, \sigma_{10}, \sigma_{11}) & T_1 & T_2
  & \text{rank } 1 & \text{rank } 2 \\
  \hline
p < -1     & (2, -2, 2, -2) & 0 & 0 & \text{no}  & \text{no}  \\
-1 < p < 0 & (2, -2, 0, 0)  & 1 & 1 & \text{no}  & \text{no}  \\
p = 0      & (2, -1, 0, -1) & 1 & 1 & \text{no}  & \text{no}  \\
0 < p < 1  & (2, 0, 0, -2)  & 1 & 1 & \text{no}  & \text{no}  \\
1 < p < 2  & (2, 0, 0, 2)   & 1 & 2 & \text{yes} & \text{yes} \\
p = 2      & (2, 1, 0, 1)   & 1 & 2 & \text{yes} & \text{yes} \\
2 < p < 3  & (2, 2, 0, 0)   & 1 & 2 & \text{yes} & \text{yes} \\
p > 3      & (2, 2, 2, 2)   & 0 & 2 & \text{no}  & \text{yes} \\
  \hline
\end{array}
\]
Hence, for generic $p$, we have that $G$ is reachable from $F$ via a rank-1 update if and only if $1 < p < 3$ and reachable via a rank-2 update if and only if $p > 1$ (excluding $p = 3$).

To check this result, we can use Weyl's inequality with the eigenvalues $\beta = (p-1, p+1)$ of $G$. Since $(F,G)$ is generic, these inequalities become strict.
Plugging this into (\ref{eq:interlacing}) with $r=1$  gives
\[
0 < p-1 < 2 < p+1,
\]
or equivalently $1 < p < 3$, which matches what we computed above using Corollary \ref{cor:rank}.
With $r=2$, we get
\[
  0 < p-1 \,\, \wedge \,\, 2 < p+1,
  \]
  or equivalently $p > 1$, which also matches what we computed above.
 \end{example}

 \begin{remark}
     The values $p = \pm 1$ and $p=3$ are excluded from the above because at these values $(F,G)$ is not generic.
     However, $G$ is still spectrally reachable from $F$ at $p = 3$ (e.g. $E = \diag(4,0)$) and $p=1$ ($E = 0$), so the restriction to genericity is a limitation of Corollary \ref{cor:rank} and not of the underlying theory.
     At these $p$ values, computing $T_1$ and $T_2$ may not yield an integer; for example $T_1 = \frac12$ if $p=3$.
     Extending Corollary \ref{cor:rank} to non-generic pairs is left for future work.
 \end{remark}
   }

{  
  \section{Conclusion}
  In this section, we will summarize the contributions of this paper and discuss some potential future research directions.

  \paragraph{Summary:}
  In the present work we have tackled the following problem: 
  given parametric  real symmetric matrices $F$~and~$G$ and an eigenvalue configuration $c$, generate a condition on the parameters such that~$\EC(F,G) = c$.
  To do this, we developed an algorithm and a robust theory of eigenvalue configurations, including a natural definition for real symmetric matrices which may share eigenvalues and/or have eigenvalues with multiplicity.

  \paragraph{Future directions:}
  We are investigating ways to prune and/or simplify the output.
  The output is essentially a disjunction of conjunctions of polynomial equalities or inequalities in the parameters. 
  When the input matrices depend on many parameters, e.g. when all entries are independent parameters, it is likely that every conjunction is consistent; that is, there exist values for the parameters that make the conjunction true.
  However, when the input matrices depend on only a few parameters, then it is likely that many of the conjunctions will be inconsistent; that is, no choice of values for the parameters can make these conjunctions true.
  Many real-world problems involve inputs with only a few parameters.
  Thus, we would like to study ways of systematically pruning inconsistent conjunctions from the output.


\bigskip\noindent\textbf{Acknowledgements.} Hoon Hong was partially supported
by US National Science Foundation NSF-CCF-2212461 and CCF 2331401. {Generative AI was used in some sections to aid in formatting and structuring, and was used to aid idea generation for Section \ref{sec:update}.}

\bibliographystyle{unsrt}
\bibliography{../../../reference/refs}

\end{document}